\documentclass[reqno, 11pt]{amsart}

\usepackage{amsmath,amssymb,amsthm,amsfonts,mathrsfs}
\usepackage{fullpage}
\usepackage{xcolor}
\usepackage{graphicx}
\usepackage{makecell}
\usepackage{listings}
\usepackage{bbm}

\definecolor{codegreen}{rgb}{0,0.6,0}
\definecolor{codegray}{rgb}{0.5,0.5,0.5}
\definecolor{codepurple}{rgb}{0.58,0,0.82}
\definecolor{backcolour}{rgb}{0.95,0.95,0.92}

\lstdefinestyle{mystyle}{
  backgroundcolor=\color{backcolour},   commentstyle=\color{codegreen},
  keywordstyle=\color{magenta},
  numberstyle=\tiny\color{codegray},
  stringstyle=\color{codepurple},
  basicstyle=\ttfamily\footnotesize,
  breakatwhitespace=false,         
  breaklines=true,                 
  captionpos=b,                    
  keepspaces=true,                 
  numbers=left,                    
  numbersep=5pt,                  
  showspaces=false,                
  showstringspaces=false,
  showtabs=false,                  
  tabsize=2
}

\usepackage{listings}
\usepackage[colorlinks=true,
linkcolor=violet,
anchorcolor=blue,
citecolor=violet
]{hyperref}
\allowdisplaybreaks % allow equations break between two pages
\newtheorem{theorem}{Theorem}[section]
\newtheorem{lemma}[theorem]{Lemma}
\newtheorem{corollary}[theorem]{Corollary}
\newtheorem{proposition}[theorem]{Proposition}

\newtheorem*{conjecture*}{Conjecture}
\newtheorem*{theorem*}{Theorem}
\theoremstyle{definition}

\theoremstyle{remark}

\newtheorem*{remark*}{Remark}
\newcommand{\Z}{\mathbb{Z}}

\newcommand{\R}{\mathbb{R}}
\newcommand{\C}{\mathbb{C}}

\newcommand{\cE}{\mathcal{E}}
\newcommand{\cZ}{\mathcal{Z}}
\newcommand{\cN}{\mathcal{N}}
\newcommand{\cS}{\mathcal{S}}
\newcommand{\cQ}{\mathcal{Q}}
\newcommand{\cJ}{\mathcal{J}}

\newcommand{\ov}[1]{\overline{#1}}
\newcommand{\sC}{\mathsf{C}}
\newcommand{\sD}{\mathsf{D}}
\newcommand{\sN}{\mathsf{N}}
\newcommand{\sR}{\mathsf{R}}
\newcommand{\sS}{\mathsf{S}}
\newcommand{\sq}{\mathsf{q}}
\newcommand{\sK}{\mathsf{K}}
\DeclareMathOperator{\sm}{sm}
\newcommand{\dd}{\mathrm{d}}

\renewcommand{\le}{\leqslant}

\usepackage{colonequals}

\author{\textbf{Runbo Li}}
\address{Department of Mathematics, Imperial College London, London SW7 2AZ, UK}
\email{runbo.li26@imperial.ac.uk, runbo.li.carey@gmail.com}

\makeatletter
\@namedef{subjclassname@2020}{\textup{}2020 Mathematics Subject Classification}
\makeatother

\title[]{Primes in arithmetic progressions to large moduli with multilinear weights}
\subjclass[2020]{\textbf{11N05}, \textbf{11N13}, \textbf{11N36}, \textbf{11L05}}
\keywords{\textbf{Primes in arithmetic progressions}, \textbf{Bombieri--Vinogradov theorem}, \textbf{dispersion method}, \textbf{sums of Kloosterman sums}}

\begin{document}

\begin{abstract}
We study the average distribution of primes up to $x$ in arithmetic progressions to moduli $q$ larger than $\sqrt{x}$, where $q = q_1 \cdots q_n$ ($n \geqslant 2$) with conveniently sized divisors $q_1, \ldots, q_n$ and the average is weighted by an arbitrary divisor-bounded multilinear weight. Specifically, we show that the primes up to $x$ are equidistributed in the residue class $a \pmod{q_1 q_2}$ on average over moduli $q_1 q_2$ as large as $x^{8/15 - \varepsilon}$ when the average is weighted by a bilinear weight and the sizes of $q_1$ and $q_2$ lie in an explicit region, improving the exponent $9/17$ for bilinear weights obtained recently by Li. The new ingredient here is an estimate for the Type-II convolution sums with both factors have size about $x^{1/2}$. With this estimate, the level of distribution $8/15$ in the bilinear case is no longer limited by the Type-II sums but by the Type-I$_3$ sums, or the current range in the ternary divisor problem.
We also get further improvements for $n \geqslant 3$ using a simple form of Harman's sieve together with various arithmetic information, including the Type-II estimate mentioned above. We group the divisors $q_1, \ldots, q_n$ in different ways in order to use arithmetic information more efficiently. Specifically, we reach $x^{145/266 - \varepsilon}$ for trilinear weights, $x^{457/814 - \varepsilon}$ for quadrilinear weights and $x^{3614/6353 - \varepsilon}$ for quintilinear weights. As the number of divisors increases, the exponent approaches $4/7$, the long-standing result with well-factorable weights proved by Bombieri, Friedlander and Iwaniec.
\end{abstract}

\maketitle

%\tableofcontents

\section{Introduction}

\subsection{Background}
For coprime integers $a$ and $q > 1$, let $\pi(x; q, a)$ denote the number of primes $p \leqslant x$ with $p \equiv a \pmod q$. The famous Bombieri--Vinogradov theorem \cite{BombieriBV} \cite{VinogradovBV} states that for every $A > 0$ there is $B = B(A)$ such that
\begin{equation}\label{eq:BV}
\sum_{q \leqslant x^{1/2} (\log x)^{-B} } \max_{(a, q) = 1} \left| \pi(x; q, a) - \frac{\pi(x)}{\varphi(q)} \right| \ll \frac{x}{(\log x)^A}.
\end{equation}
The Elliott--Halberstam conjecture \cite{ElliottHalberstam} predicts that the range of summation may be extended to $q \leqslant x^{1 - \varepsilon}$, but no estimate of the strength of \eqref{eq:BV} is known with $x^{1/2}$ replaced by $x^{1/2 + \delta}$ for any $\delta > 0$.

Beginning with the work of Fouvry and Iwaniec \cite{FouvryIwaniec} \cite{FouvryA1} and of Bombieri, Friedlander and Iwaniec \cite{BFI} \cite{BFI2} \cite{BFI3}, one has nevertheless been able to go beyond $x^{1/2}$ for a \textit{fixed} non-zero residue class $a$, provided that the absolute values in \eqref{eq:BV} are replaced by weights with suitable factorisation properties, or that the moduli are restricted. Thus one knows such results for well-factorable weights of level $x^{4/7 - \varepsilon}$ \cite{BFI}, for triply well-factorable weights of level $x^{3/5 - \varepsilon}$ \cite{MaynardLargeModuliII} and $x^{5/8 - \varepsilon}$ \cite{Pascadi2}, for smooth moduli \cite{ZhangYitang} \cite{Polymath8a} \cite{Stadlmann525}, and for moduli with a conveniently sized divisor: Maynard \cite{MaynardLargeModuliI} proved that
\begin{equation}\label{eq:Maynard}
\sum_{\substack{q_1 \leqslant Q_1 \\ q_2 \leqslant Q_2 \\ (q_1 q_2, a) = 1}} \left| \pi(x; q_1 q_2, a) - \frac{\pi(x)}{\varphi(q_1 q_2)} \right| \ll \frac{x}{(\log x)^A}
\end{equation}
whenever 
$$
Q_1^2 Q_2 < x^{1 - \varepsilon}, \quad Q_1^7 Q_2^{12} < x^{4 - \varepsilon}, \quad Q_1^{19} Q_2^{20} < x^{10 - \varepsilon},
$$
which allows $Q_1 Q_2$ to be as large as $x^{11/21 - \varepsilon}$.

An intermediate situation between \eqref{eq:Maynard} and well-factorable weights is that of \textit{bilinear weights}: one asks for estimates of the form
\begin{equation}\label{eq:bilinear}
\sum_{\substack{q_1 \leqslant Q_1 \\ q_2 \leqslant Q_2 \\ (q_1 q_2, a) = 1}} \lambda_1(q_1) \lambda_2(q_2) \left( \pi(x; q_1 q_2, a) - \frac{\pi(x)}{\varphi(q_1 q_2)} \right) \ll \frac{x}{(\log x)^A}
\end{equation}
for arbitrary divisor-bounded sequences $\lambda_1, \lambda_2$. In 1986, Bombieri, Friedlander and Iwaniec \cite{BFI} first showed \eqref{eq:bilinear} holds if
$$
Q_1 < x^{1/3}, \quad Q_2 < x^{1/5}, \quad Q_1^{5} Q_2^{2} < x^{2}, \quad Q_1 Q_2 < x^{29/56}.
$$
In 1987, Fouvry \cite{FouvryA2} proved that \eqref{eq:bilinear} holds if either
$$
Q_1 Q_2^{3} < x, \quad Q_1 Q_2 < x^{29/56}, \quad Q_1 < \max\left(x^{1/2} Q_2^{-1}, x^{2/5} Q_2^{-2/5} \right)
$$
or
$$
Q_1 Q_2^{3} < x, \quad Q_1 Q_2 < x^{29/56}, \quad Q_1^{4} Q_2 < x^{403/266}, \quad Q_1^{7/4} Q_2 < x^{403/532}.
$$
In 1998, Baker and Harman \cite{677} improved the result of Bombieri, Friedlander and Iwaniec \cite{BFI} by removing the condition $Q_1^{5} Q_2^{2} < x^{2}$ in their result above. They showed that \eqref{eq:bilinear} holds if
$$
Q_1 < x^{1/3}, \quad Q_2 < x^{1/5}, \quad Q_1 Q_2 < x^{29/56}.
$$
Clearly \eqref{eq:Maynard} implies \eqref{eq:bilinear} in Maynard's region. Recently Li \cite{LRB2602} showed, by combining the arithmetic information of \cite{MaynardLargeModuliI}, \cite{FouvryA2}, \cite{IwaniecPomykala} and \cite{Lichtman2} with refinements of Harman's sieve, that \eqref{eq:bilinear} holds whenever
\begin{equation}\label{eq:Li}
Q_1^2 Q_2 < x^{1 - \varepsilon}, \qquad Q_1^7 Q_2^{12} < x^{4 - \varepsilon},
\end{equation}
so that $Q_1 Q_2$ may be as large as $x^{9/17 - \varepsilon}$. For quadrilinear weights Lichtman \cite{Lichtman2} reached moduli of size $x^{17/32 - \varepsilon}$, and Li \cite{LRB2602} obtained the same exponent $17/32$ for trilinear weights.

\subsection{Main Results}
One of the main purposes of this paper is to show that for bilinear weights one can reach the exponent $8/15$.
\begin{theorem}\label{thm:main}
Let $a$ be a fixed non-zero integer and let $\varepsilon > 0$, $A > 0$, $B_0 > 0$. Let $Q_1, Q_2 \geqslant 1$ satisfy
\begin{equation}\label{eq:region}
Q_1^2 Q_2 \leqslant x^{1 - \varepsilon}, \qquad Q_1^8 Q_2^{14} \leqslant x^{5 - \varepsilon}, \qquad Q_1^{16} Q_2^{107} \leqslant x^{41/2 - \varepsilon},\qquad Q_1 Q_2 \leqslant x^{8/15 - \varepsilon}.
\end{equation}
Let $\lambda_1, \lambda_2$ be complex sequences with $|\lambda_i(q)| \leqslant \tau(q)^{B_0}$. Then
$$
\sum_{\substack{q_1 \leqslant Q_1 \\ q_2 \leqslant Q_2 \\ (q_1 q_2, a) = 1}} \lambda_1(q_1) \lambda_2(q_2) \left( \pi(x; q_1 q_2, a) - \frac{\pi(x)}{\varphi(q_1 q_2)} \right) \ll_{a, \varepsilon, A, B_0} \frac{x}{(\log x)^A}.
$$
\end{theorem}
Write $Q_i = x^{\theta_i}$. In terms of $\theta_1, \theta_2$, the region \eqref{eq:region} is given by $2 \theta_1 + \theta_2 < 1$, $8 \theta_1 + 14 \theta_2 < 5$, $16 \theta_1 + 107 \theta_2 < 41/2$ and $\theta_1 + \theta_2 < 8/15$. The first condition forces $\theta_1 < 1/2$, so that $q_1$ is the ``large'' and $q_2$ the ``small'' factor whenever $\theta_1 + \theta_2 \geqslant 1/2$. Parametrising by the level, we obtain the following.
\begin{corollary}\label{cor:level}
Let $0 < \delta < 1/30$ and let $\theta_2$ satisfy
$$
2 \delta < \theta_2 < \min\left(\frac{1 - 8 \delta}{6},\ \frac{25 - 32 \delta}{182}\right).
$$
Then \eqref{eq:bilinear} holds with $Q_1 = x^{1/2 + \delta - \theta_2}$ and $Q_2 = x^{\theta_2}$. In particular, for every fixed $\theta_2 \in (1/15, 11/90)$ and every sufficiently small $\varepsilon > 0$, the estimate \eqref{eq:bilinear} holds with $Q_2=x^{\theta_2}$ and $Q_1 Q_2=x^{8/15 - \varepsilon}$.
\end{corollary}
Indeed, with $\theta_1 = 1/2 + \delta - \theta_2$, the four conditions read $\theta_2 > 2 \delta$, $\theta_2 < (1 - 8 \delta)/6$, $\theta_2 < (25 - 32 \delta)/182$ and $\delta < 1/30$, and Theorem~\ref{thm:main} applies with a sufficiently small $\varepsilon$. A simple verification shows that Li's $9/17$ region \eqref{eq:Li} with $Q_1 Q_2 > x^{1/2 - \varepsilon}$ (the part with $Q_1 Q_2 \leqslant x^{1/2 - \varepsilon}$ being covered by the Bombieri--Vinogradov theorem) is contained in \eqref{eq:region}: if $2 \theta_1 + \theta_2 < 1$, $7 \theta_1 + 12 \theta_2 < 4$ and $\theta_1 + \theta_2 \geqslant 1/2$, then $\theta_2 < 1/10$ and $\theta_1 > 2/5$, whence $8 \theta_1 + 14 \theta_2 < 72/17 + 3/5 = 411/85 < 5$, $16 \theta_1 + 107 \theta_2 < 144/17 + 91/10 = 2987/170 < 41/2$ and $\theta_1 + \theta_2 < 9/17 < 8/15$. The known results are compared in Table~\ref{tab:comparison}.

\begin{table}[ht]
\renewcommand{\arraystretch}{1.25}
\begin{tabular}{llll}
\hline
& Type of average & Conditions on $\theta_i$ & $\sup(\sum \theta_i)$ \\
\hline
B--F--I \cite{BFI} & bilinear weights & $\theta_1 < \frac{1}{3}$, $\theta_2 < \frac{1}{5}$, & $29/56 \approx 0.5179$ \\
& & $5 \theta_1 + 2 \theta_2 < 2$, $\theta_1 + \theta_2 < \frac{29}{56}$ & \\
Fouvry \cite{FouvryA2} & bilinear weights & $\theta_1 + 3 \theta_2 < 1$, $\theta_1 + \theta_2 < \frac{29}{56}$, & $29/56 \approx 0.5179$ \\
& & $\theta_1 < \max(\frac{1}{2} - \theta_2, \frac{2}{5}(1 - \theta_2))$; & \\
& & $\theta_1 + 3 \theta_2 < 1$, $\theta_1 + \theta_2 < \frac{29}{56}$, & $29/56 \approx 0.5179$ \\
& & $4 \theta_1 + \theta_2 < \frac{403}{266}$, $\frac{7}{4} \theta_1 + \theta_2 < \frac{403}{532}$ & \\
B--H \cite{677} & bilinear weights & $\theta_1 < \frac{1}{3}$, $\theta_2 < \frac{1}{5}$, $\theta_1 + \theta_2 < \frac{29}{56}$ & $29/56 \approx 0.5179$ \\
Maynard \cite{MaynardLargeModuliI} & absolute values & $2 \theta_1 + \theta_2 < 1$, $7 \theta_1 + 12 \theta_2 < 4$, & $11/21 \approx 0.5238$ \\
& & $19 \theta_1 + 20 \theta_2 < 10$ & \\
Lichtman \cite{Lichtman2} & quadrilinear weights & $\theta_1 + \theta_2 < \frac{1}{2}$, $\theta_1 + 2 \theta_3 < \frac{1}{2}$, & $17/32=0.53125$ \\
& & $2 \theta_3 < \theta_2 < \frac{1}{32}$ & \\
Li \cite{LRB2602} & bilinear weights & $2 \theta_1 + \theta_2 < 1$, $7 \theta_1 + 12 \theta_2 < 4$ & $9/17 \approx 0.5294$ \\
& trilinear weights & complicated, see \cite{LRB2602} & $17/32=0.53125$ \\
\hline
Theorem~\ref{thm:main} & bilinear weights & $2 \theta_1 + \theta_2 < 1$, $8 \theta_1 + 14 \theta_2 < 5$, & $8/15 \approx 0.5333$ \\
& & $16 \theta_1 + 107 \theta_2 < \frac{41}{2}$, $\theta_1 + \theta_2 < \frac{8}{15}$ & \\
Theorems~\ref{thm:trilinear}--\ref{thm:345linear} & trilinear weights & complicated & $145/266 \approx 0.5451$ \\
Theorem~\ref{thm:345linear} & quadrilinear weights & complicated & $457/814 \approx 0.5614$ \\
Theorem~\ref{thm:345linear} & quintilinear weights & complicated & $3614/6353 \approx 0.5688$ \\
\hline
B--F--I \cite{BFI} & well-factorable weights & & $4/7 \approx 0.5714$ \\
\hline
\end{tabular}
\medskip
\caption{\small{Average equidistribution results for primes in a fixed residue class to factored moduli}}\label{tab:comparison}
\end{table}

Since a Dirichlet convolution of divisor-bounded sequences is also divisor-bounded, Theorem~\ref{thm:main} immediately covers multilinear weights.
\begin{corollary}\label{cor:multilinear}
Let $k \geqslant 2$, let $Q_1, \ldots, Q_k \geqslant 1$ and let $\lambda_1, \ldots, \lambda_k$ be complex sequences with $|\lambda_i(q)| \leqslant \tau(q)^{B_0}$. Suppose that for some $I \subset \{1, \ldots, k \}$, the pair $\left(\prod_{i \in I} Q_i, \prod_{i \notin I} Q_i \right)$ satisfies \eqref{eq:region}. Then, for every fixed non-zero integer $a$ and every $A > 0$,
$$
\sum_{\substack{q_i \leqslant Q_i\ (1 \leqslant i \leqslant k) \\ (q_1 \cdots q_k, a) = 1}} \lambda_1(q_1) \cdots \lambda_k(q_k) \left(\pi(x; q_1 \cdots q_k, a) - \frac{\pi(x)}{\varphi(q_1 \cdots q_k)} \right) \ll \frac{x}{(\log x)^A}.
$$
\end{corollary}
\begin{proof}
Put $\lambda^I(q) := \sum \prod_{i \in I} \lambda_i(q_i)$, the sum being over all factorizations $q = \prod_{i \in I} q_i$ with $q_i \leqslant Q_i$, and define $\lambda^J$ similarly with respect to the complement of $I$. Then $|\lambda^I(q)|, |\lambda^J(q)| \leqslant \tau(q)^{k(B_0+1)}$, the sequences are supported on $q \leqslant \prod_{i \in I} Q_i$ and $q \leqslant \prod_{i \notin I} Q_i$ respectively, and the sum in question equals the sum of Theorem~\ref{thm:main} with $(\lambda_1, \lambda_2)$ replaced by $(\lambda^I, \lambda^J)$.
\end{proof}
Indeed, we can obtain better results for multilinear weights by grouping divisors in different ways and applying Harman's sieve, and the details will be given in Section~\ref{sec:multilinear}.

\subsection{A new Type-II estimate near \texorpdfstring{$x^{1/2}$}{the square root of x}}
Results of the above type are proved by decomposing the indicator function of the primes into convolutions, and the range of moduli is limited by the available estimates for
\begin{equation}\label{eq:Delta}
\Delta(q) = \Delta_{\alpha, \beta}(q; a) := \sum_{m \sim M} \sum_{n \sim N} \alpha_n \beta_m \left(\mathbbm{1}_{m n \equiv a (\operatorname{mod} q)}-\frac{\mathbbm{1}_{(m n, q) = 1}}{\varphi(q)} \right), \qquad M N \asymp x,
\end{equation}
in various ranges of $N$. When the moduli have size $x^{1/2 + \delta}$, the most delicate range is $x^{1/2 - \delta} \leqslant N \leqslant x^{1/2 + \delta}$, where no estimate is known for general moduli. In the works of Zhang \cite{ZhangYitang}, Polymath8a \cite{Polymath8a} and Maynard \cite{MaynardLargeModuliI} this range is treated with the help of the Weil bound for (incomplete) Kloosterman sums, and for moduli $q_1 q_2$ with $q_1 \sim Q$, $q_2 \sim R$ this requires $Q^7 R^{12} < x^{4}$ (see [\cite{MaynardLargeModuliI}, Proposition 8.2]). This is the origin of the second condition in \eqref{eq:Li}, which is the reason why the region \eqref{eq:Li} stops at $9/17$. Our main new result removes this obstruction for bilinear weights.

Throughout the paper we put
$$
z_0 := x^{1/(\log\log x)^3}, \qquad y_0 := x^{1/\log\log x},
$$
and we say that a complex sequence $(\alpha_n)$ supported on $n \sim N$ satisfies the \textit{Siegel--Walfisz condition} if for every $d \geqslant 1$, $k \geqslant 1$, $(b, k) = 1$ and every $A > 1$ we have
\begin{equation}\label{eq:SW}
\left|\sum_{\substack{n \sim N \\ n \equiv b (\operatorname{mod} k) \\ (n, d) = 1 }} \alpha_n - \frac{1}{\varphi(k)} \sum_{\substack{n \sim N \\ (n, d k) = 1}} \alpha_n \right| \ll_A \frac{N \tau(d)^{B_0}}{(\log N)^A}.
\end{equation}
We note that $\alpha_n$ certainly satisfies the Siegel--Walfisz condition \eqref{eq:SW} if $\alpha_n = 1$, if $\alpha_n = \mu(n)$, or if
$$
\alpha_n = \sum_{\substack{p_1 \cdots p_j = n \\ p_i \sim P_i}} 1
$$
by the Siegel--Walfisz theorem.

For $q \in \mathbb{N}$, let $\mathit{\lambda_1} = \mathit{\lambda_1}(q)$ be the smallest positive eigenvalue of the hyperbolic Laplacian for the congruence subgroup $\Gamma_0(q)$, and put $\vartheta_q := \max(0, \sqrt{1 - 4 \mathit{\lambda_1}})$. We say that $\vartheta \in [0, 1/2]$ is \textit{admissible} if $\vartheta_q \leqslant \vartheta$ for all $q \in \mathbb{N}$. Selberg's eigenvalue conjecture \cite{SelbergEigenvalue} asserts that $\vartheta = 0$ is admissible, while his bound $\mathit{\lambda_1} \geqslant 3/16$ gives $\vartheta = 1/2$ is admissible. The best unconditional result now is due to Kim and Sarnak \cite{KimSarnak}, where they proved that $\vartheta = 7/32$ is admissible.

\begin{theorem}[Type-II estimate near $x^{1/2}$]\label{thm:typeII}
Let $\vartheta$ be admissible. Let $a$ be a fixed non-zero integer and let $\varepsilon > 0$, $A > 0$, $B_0 > 0$. Let $M, N, Q, R \geqslant 1$ satisfy $M N \asymp x$ and
\begin{align}
& Q x^{\varepsilon} \leqslant N \leqslant x^{1/2 + \varepsilon/100}, \label{eq:T1} \\
& N^{1/2} Q^{1/2} R^{4} \leqslant x^{1 - \varepsilon}, \label{eq:T2} \\
& N^{1/2 + \vartheta} Q^{1/2} R^{4 - 3 \vartheta} \leqslant x^{1 - \varepsilon}. \label{eq:T3}
\end{align}
Let $\alpha_n, \beta_m, \lambda_q, \nu_r$ be complex sequences with $|\alpha_n|, |\beta_n|, |\lambda_n|, |\nu_n| \leqslant \tau(n)^{B_0}$, and assume that $\alpha_n$ is supported on integers $n$ with $P^{-}(n) \geqslant z_0$ and satisfies the Siegel--Walfisz condition \eqref{eq:SW}. Then, with $\Delta(q)$ as in \eqref{eq:Delta},
$$
\sum_{q \sim Q} \sum_{\substack{r \sim R \\ (q r,a) = 1}} \lambda_q \nu_r \Delta(q r) \ll_{a, \varepsilon, A, B_0} \frac{x}{(\log x)^A}.
$$
\end{theorem}
With $\vartheta = 7/32$, the condition \eqref{eq:T3} reads
$$
N^{23/32} Q^{1/2} R^{107/32} \leqslant x^{1 - \varepsilon}.
$$
When both $\alpha$ and $\beta$ are of the above type, the roles of $m$ and $n$ may be interchanged, and we obtain an estimate in a range of $N$ which is symmetric about $x^{1/2}$.
\begin{corollary}\label{cor:symmetric}
Let $\varepsilon > 0$ and let $Q, R \geqslant 1$ satisfy
$$
Q R^8 \leqslant x^{3/2 - 3 \varepsilon}, \qquad Q^{16} R^{107} \leqslant x^{41/2 - 17 \varepsilon}.
$$
Let $M N \asymp x$ with $Q x^{\varepsilon} \leqslant N \leqslant x^{1 - \varepsilon}/Q$, and let $\alpha_n, \beta_m, \lambda_q, \nu_r$ be divisor-bounded as in Theorem~\ref{thm:typeII}, where now both $\alpha_n$ and $\beta_m$ are supported on integers without prime factors below $z_0$ and satisfy the Siegel--Walfisz condition \eqref{eq:SW}. Then
$$
\sum_{q \sim Q} \sum_{\substack{r \sim R \\ (q r,a) = 1}} \lambda_q \nu_r \Delta(q r) \ll \frac{x}{(\log x)^A}.
$$
\end{corollary}

\subsection{Outline of the method}\label{subsec:outline}
We describe the proof of Theorem~\ref{thm:typeII} in a simplified situation, ignoring coprimality conditions and factors $x^{o(1)}$. Following Linnik's dispersion method, we apply Cauchy's inequality to the sum over $q$ and $m$; this removes the unknown coefficients $\lambda_q$ and $\beta_m$ and makes $q$ and $m$ \textit{smooth} variables. After an application of Poisson summation to the sum over $m$, and after the removal of the expected main terms, one has to show that
\begin{equation}\label{eq:model}
\sum_{q \sim Q} \sum_{r_1, r_2 \sim R} \frac{\nu_{r_1} \ov{\nu_{r_2}}}{q r_1 r_2} \sum_{\substack{n_1, n_2 \sim N \\ n_1 \equiv n_2 (\operatorname{mod} q) }} \alpha_{n_1} \ov{\alpha_{n_2}} \sum_{1 \leqslant |h| \leqslant H} e\left(a h \frac{n_1 - n_2}{q} \cdot \frac{\ov{n_2 r_1}}{n_1 r_2} \right), \qquad H \approx \frac{Q R^2}{M},
\end{equation}
is smaller than the trivial bound $N^2 H/Q$ by a factor a little larger than $H$.

The treatments of sums like \eqref{eq:model} in the literature differ in the choice of the variables which are kept outside in a second application of Cauchy's inequality. In the first method of [\cite{BFI}, Section 8] one replaces the smooth variable $q$ by the complementary divisor $f=(n_1 - n_2)/q$ of $n_1 - n_2$, and keeps $(n_1 r_2, n_2, f)$ outside and $(h, r_1)$ inside; the off-diagonal terms are sums of Kloosterman fractions over the two smooth variables $c = n_1 r_2$ and $d = n_2$ and are estimated by the bounds of Deshouillers and Iwaniec \cite{DI82Kl}, but the diagonal terms are admissible only if $N^2 R < x$. In the work of Zhang \cite{ZhangYitang} (see [\cite{MaynardLargeModuliI}, Chapter 17]) the variables $(q, n_1, n_2)$ are kept outside and $(h, r_1, r_2)$ inside; the diagonal terms are then admissible for all $N$ up to about $x^{1/2}$, but the coefficients depend on $q$ (as they must in a result with absolute values), and the off-diagonal terms are estimated by completing the sums over $n_1, n_2$ and applying the Weil bound, which leads to $Q^7 R^{12} < x^4$.

Our observation is that for bilinear weights the two approaches can be combined. Since $q$ is smooth we may pass to $f = (n_1-n_2)/q$, which is a \textit{short} variable, $|f| \leqslant N/Q$. We then keep $(f, n_1, n_2)$ outside and $(h, r_1, r_2)$ inside, and we drop by positivity the congruence $n_1 \equiv n_2 \pmod{f}$, so that $n_1$ and $n_2$ become free smooth variables of length $N$. After expanding the square, and after the usual decomposition of the pairs $(r_1, r_2)$ into classes within which all coprimality conditions hold \cite{FouvryA1}, the off-diagonal terms are
$$
\sum_{f} \sum_{\substack{(h, r_1, r_2) \\ (h', r_1', r_2') }} \xi(h, r_1, r_2) \ov{\xi(h', r_1', r_2')} \sum_{n_1, n_2} \Psi\left(\frac{n_1}{N}\right) \Psi\left(\frac{n_2}{N}\right) e\left(a f k \frac{\ov{n_2 r_1 r_1'}}{n_1 [r_2, r_2']}\right)
$$
with certain bounded coefficients $\xi$ and with
$$
k = \frac{h r_1' r_2' - h' r_1 r_2}{(r_2, r_2')} \neq 0.
$$
This is exactly a sum of the type considered in [\cite{DI82Kl}, Theorem 12], with both ``level'' variables $r = r_1 r_1'$ and $s = [r_2, r_2']$ present. The resulting bound is admissible under \eqref{eq:T2} and \eqref{eq:T3}; these two conditions come from the regular and from the (possible) exceptional spectrum respectively. The diagonal terms $k = 0$ are admissible as long as $N \leqslant x^{1/2-\varepsilon}$. In order to reach $N = x^{1/2+\varepsilon/100}$ we use the de-amplification device of Maynard [\cite{MaynardLargeModuliI}, Chapter 14]: before the first application of Cauchy's inequality, the variable $n$ is split into residue classes modulo $e$ for auxiliary moduli $e \sim E = x^{O(\varepsilon)}$, which produces the stronger congruence $n_1 \equiv n_2 \pmod{q e}$ in \eqref{eq:model}. The congruence modulo $e$ cannot be dropped; we retain it in the outer summation and appeal to the generalization of the Deshouillers--Iwaniec bound to sums with congruence conditions due to Lichtman \cite{Lichtman3}, where the dependence on the parameter $\vartheta$ is explicit.

\subsection{Scope and limitations}\label{subsec:scope}
With Theorem~\ref{thm:typeII} in hand the Type-II sums are under control for moduli of size up to about $x^{223/396}$ (if $Q^2 R \approx x$). The exponent $8/15$ in Theorem~\ref{thm:main} comes from a different source: products of three primes of size about $x^{1/3}$ are handled by means of estimates for the ternary divisor function [\cite{MaynardLargeModuliI}, Chapter 20], which are currently available only for moduli up to $x^{8/15}$. The conditions $Q_1^8 Q_2^{14} \leqslant x^{5 - \varepsilon}$ and $Q_1^{16} Q_2^{107} \leqslant x^{41/2 - \varepsilon}$ are needed when applying other arithmetic information (the second one is also needed in Theorem~\ref{thm:typeII}), and $Q_1^2 Q_2 \leqslant x^{1 - \varepsilon}$ guarantees that the range of Theorem~\ref{thm:typeII} covers the critical interval around $x^{1/2}$.

The method relies in an essential way on the fact that the coefficient of the large modulus is removed by the first application of Cauchy's inequality. It therefore gives nothing new for general moduli (the exponent $1/2$ in \eqref{eq:BV}, for either fixed or uniform residue classes), nor for estimates with absolute values such as \eqref{eq:Maynard}; the region of \cite{BFI} \cite{677}, where $q_1$ and $q_2$ are of comparable size, is also disjoint from \eqref{eq:region}, and $29/56$ remains the best exponent there.

Sections \ref{sec:prelim}--\ref{sec:proofTypeII} contain the proof of Theorem~\ref{thm:typeII}. The new arguments are in Section~\ref{sec:E2}; the reduction to exponential sums in Section~\ref{sec:dispersion} is a variant of [\cite{MaynardLargeModuliI}, Proposition 14.4], and we indicate precisely where it differs. In Section~\ref{sec:main} we deduce Theorem~\ref{thm:main}, following the architecture of \cite{MaynardLargeModuliI}. In Section~\ref{sec:multilinear} we further apply a simple form of Harman's sieve to get better results for multilinear weights.
The only new arithmetic information input is Theorem~\ref{thm:typeII}, the remaining ingredients are taken from \cite{BFI}, \cite{BFI2}, \cite{FouvryA1}, \cite{FouvryA2}, \cite{IwaniecPomykala}, \cite{Lichtman2}, \cite{Lichtman3}, \cite{MaynardLargeModuliI} and \cite{Pascadi2}.

\subsection*{Notation}
Throughout this paper, we always suppose that $x$ is a sufficiently large real number, $a$ is a fixed non-zero integer and $\varepsilon$, $\varepsilon_0$ are sufficiently small positive real numbers. We use both Landau and Vinogradov notations, and implied constants may depend on $a$, $\varepsilon$, and on the exponents $A, B_0$. We write $e(t) := \exp(2 \pi i t)$, and $n \sim N$ for $N < n \leqslant 2N$. As usual, $\tau$, $\varphi$, $\mu$ denote the divisor function, Euler's totient function and the M\"obius function, $P^{-}(n)$ and $P^{+}(n)$ denote the least and the greatest prime factor of $n$, and $(m,n)$, $[m,n]$ denote the greatest common divisor and the least common multiple. If $(n,c)=1$, then $\ov{n}$ denotes an inverse of $n$ modulo $c$; for complex coefficients a bar denotes complex conjugation. A sequence is called \textit{$1$-bounded} if all its terms have modulus at most $1$. We write $\sm(n; z) := \prod_{p \leqslant z} p^{v_p(n)}$ for the $z$-smooth part of $n$. Following \cite{MaynardLargeModuliI}, we let $\psi_0: \R \to \R$ be the fixed smooth function in [\cite{MaynardLargeModuliI}, Chapter 6]; it is supported on $[1/2, 5/2]$ and equal to $1$ on $[1,2]$, and $\widehat{\psi}_0(\xi) = \int \psi_0(t) e(-\xi t) \dd t$ is its Fourier transform. Finally we fix a smooth function $\Psi: \R \to [0,1]$ supported on $[1/2,5/2]$ with $\Psi = 1$ on $[1,2]$.

\section{Preliminary Lemmas}\label{sec:prelim}
In this section we establish various lemmas which will be used in later sections.

\subsection{Reciprocity}
\begin{lemma}[Bezout's identity]\label{lem:bezout}
Let $(q_1, q_2)=1$. Then we have
$$
\frac{a}{q_1 q_2} \equiv \frac{a \ov{q_1}}{q_2} + \frac{a \ov{q_2}}{q_1} \pmod{1}.
$$
\end{lemma}

The next lemma will be used to compute the phase arising from Poisson summation in the dispersion method. In \cite{BFI} and \cite{MaynardLargeModuliI} the corresponding computation is carried out with Lemma~\ref{lem:bezout} under the additional assumption $(\ell r_1,r_2)=1$; it is convenient for us to avoid this assumption.
\begin{lemma}\label{lem:reciprocity}
Let $\ell, r_1, r_2, n_1, n_2 \geqslant 1$ and $b$ be integers with
$$
(r_1, r_2) = (n_1, n_2) = (n_1, r_1) = (n_2, r_2) = 1 \qquad \text{and} \qquad n_1 \equiv n_2 \pmod{\ell},
$$
and put $f := (n_1 - n_2)/\ell$. Then $(n_1 n_2, \ell) = 1$, the system of congruences
$$
m n_1 \equiv b \pmod{\ell r_1}, \qquad m n_2 \equiv b \pmod{\ell r_2}
$$
has a unique solution $m \equiv \mu \pmod{\ell r_1 r_2}$, and
$$
\frac{\mu}{\ell r_1 r_2}\equiv b f \frac{\ov{n_2 r_1}}{n_1 r_2} + \frac{b}{\ell r_1 r_2 n_1} \pmod{1}.
$$
\end{lemma}
\begin{proof}
A prime dividing $\ell$ and $n_1$ divides $n_2 = n_1 - \ell f$, which is impossible; hence $(n_1 n_2, \ell) = 1$. Consequently $(n_1, \ell r_1) = (n_2, \ell r_2) = 1$, so that two solutions of the system agree modulo $\ell r_1$ and modulo $\ell r_2$, hence modulo $[\ell r_1, \ell r_2] = \ell r_1 r_2$. Our assumptions give $(n_2 r_1, n_1 r_2)=1$. Let $X$ be an integer with $X n_2 r_1 \equiv b f \pmod{n_1 r_2}$, say $X n_2 r_1 = b f + t n_1 r_2$. Multiplying by $\ell$ and using $\ell f = n_1 - n_2$ we get
$$
\ell r_1 X n_2 + b n_2 = b n_1 + \ell t n_1 r_2.
$$
Thus $n_1$ divides $(\ell r_1 X + b)n_2$, and since $(n_1, n_2)=1$ the number $\mu:=(\ell r_1 X + b)/n_1$ is an integer. By construction $\mu n_1 = \ell r_1 X + b \equiv b\pmod{\ell r_1}$, and the displayed identity gives $\mu n_2 = b + \ell t r_2 \equiv b \pmod{\ell r_2}$. Finally
$$
\frac{\mu}{\ell r_1 r_2} = \frac{\ell r_1 X+b}{\ell r_1 r_2 n_1} = \frac{X}{n_1 r_2} + \frac{b}{\ell r_1 r_2 n_1}.
$$
\end{proof}
\begin{lemma}\label{lem:combine}
Let $n_1, n_2, r_1, r_1', r_2, r_2' \geqslant 1$ be integers with $(n_2 r_1 r_1', n_1 r_2 r_2') = 1$, and let $h, h' \in \Z$. Then
$$
h \frac{\ov{n_2 r_1}}{n_1 r_2} - h' \frac{\ov{n_2 r_1'}}{n_1 r_2'} \equiv k \frac{\ov{n_2 r_1 r_1'}}{n_1 [r_2, r_2']} \pmod{1}, \qquad k := \frac{h r_1' r_2' - h' r_1 r_2}{(r_2, r_2')}.
$$
\end{lemma}
\begin{proof}
Let $s := [r_2, r_2']$ and let $Y$ be an inverse of $n_2 r_1 r_1'$ modulo $n_1 s$. Then $r_1' Y$ is an inverse of $n_2 r_1$ modulo $n_1 r_2$ and $r_1 Y$ is an inverse of $n_2 r_1'$ modulo $n_1 r_2'$. Hence the left-hand side is congruent modulo $1$ to
$$
\frac{h r_1' Y}{n_1 r_2} - \frac{h' r_1 Y}{n_1 r_2'} = \frac{Y}{n_1 s}\left(h r_1' \frac{s}{r_2} - h' r_1 \frac{s}{r_2'}\right),
$$
and $s/r_2 = r_2'/(r_2, r_2')$, $s/r_2' = r_2/(r_2, r_2')$.
\end{proof}

\subsection{Splitting into coprime sets}
\begin{lemma}\label{lem:coprime} ([\cite{MaynardLargeModuliI}, Lemma 13.2]).
Let $\cN$ be a set of pairs $(u, v)$ of positive integers with $u, v \leqslant x$, $(u, v) = 1$ and $P^{-}(uv) \geqslant z_0$. Then there is a partition $\cN = \cN_1 \sqcup \cdots \sqcup \cN_J$ with 
$$
J \ll \exp\left(O(\log \log x)^7\right)
$$
such that $(u,v')=(u',v)=1$ whenever $(u,v)$ and $(u',v')$ lie in the same class $\cN_j$.
\end{lemma}

\subsection{Kloosterman sums}
\begin{lemma}[Incomplete Kloosterman sums]\label{lem:weil}
Let $C_0 > 0$ and let $f: \R \to \C$ be a smooth function supported on $[1/2, 5/2]$ with $\|f^{(j)}\|_{\infty} \ll_j (\log x)^{j C_0}$ for all $j \geqslant 0$. Then for all integers $b$ and $1 \leqslant c \leqslant x^{C_0}$ and all real $0 < Y \leqslant x^{C_0}$ we have
$$
\sum_{\substack{n \\ (n, c) = 1}} f\left(\frac{n}{Y}\right) e\left(\frac{b \ov{n}}{c}\right) \ll x^{o(1)} \left(c^{1/2} + \frac{Y}{c}(b,c)\right).
$$
\end{lemma}
\begin{proof}
If $Y < 1$, the sum has $O(1)$ terms and the bound is trivial. For $Y \geqslant 1$, this estimate follows from completion of the sum and the Weil bound; see [\cite{MaynardLargeModuliI}, Lemma 16.1].
\end{proof}

The following deep estimate for sums of Kloosterman sums is the main tool of this paper. It originates in the work of Deshouillers and Iwaniec [\cite{DI82Kl}, Theorem 12], where it is deduced from the Kuznetsov trace formula. The version with congruence conditions and explicit $\vartheta$-dependent is due to Lichtman \cite{Lichtman3}, with $a = 1$.
\begin{theorem}\label{thm:DI} ([\cite{Lichtman3}, Theorem 1.8]).
Let $\vartheta$ be admissible. Let $\sC, \sD, \sN, \sR, \sS \geqslant 1$, and let $\sq, c_0, d_0 \geqslant 1$ be integers with $(c_0 d_0, \sq) = 1$. Let $(b_{n,r,s})$ be a complex sequence supported on $1 \leqslant n \leqslant \sN$, $\sR < r \leqslant 2 \sR$, $\sS < s \leqslant 2 \sS$. Let $g : \R_{>0}^5 \to \C$ be a smooth function, compactly supported in $(\sC, 2\sC] \times (\sD, 2\sD] \times \R_{>0}^3$, such that
$$
\frac{\partial^{\nu_1+\cdots+\nu_5} g}{\partial c^{\nu_1} \partial d^{\nu_2} \partial n^{\nu_3} \partial r^{\nu_4} \partial s^{\nu_5}}(c,d,n,r,s) \ll_{\nu_1,\ldots,\nu_5} \left(c^{-\nu_1} d^{-\nu_2} n^{-\nu_3} r^{-\nu_4} s^{-\nu_5} \right)^{1 - \epsilon_0}
$$
for some small $\epsilon_0 > 0$ and all fixed $\nu_j \geqslant 0$. Then
$$
\sum_{\substack{c, d, n, r, s \\ c \equiv c_0, d \equiv d_0 (\operatorname{mod} \sq) \\ (\sq r d, s c) = 1}} b_{n,r,s} g(c,d,n,r,s) e\left(n \frac{\ov{rd}}{sc}\right) \ll_{\epsilon, \epsilon_0}(\sq \sC \sD \sN \sR \sS)^{\epsilon + O(\epsilon_0)} \sq^{3/2} \sK \|b\|_2,
$$
where $\|b\|_2 := \left(\sum_{n,r,s}|b_{n,r,s}|^2\right)^{1/2}$ and
\begin{equation}\label{eq:K}
\sK^2:=\sq (\sC \sS (\sN + \sR \sS) (\sC + \sD \sR) + \sN \sR \sS) + \sC \sS (\sC \sD \sqrt{\sR})^{2 \vartheta} (\sN + \sR \sS)^{1 - \vartheta} (\sC + \sD \sR)^{1 - 2 \vartheta} + \sD^2 \sN \sR.
\end{equation}
\end{theorem}
\begin{remark*}
Taking $\vartheta = 1/2$, this becomes [\cite{DrappeauKl}, Theorem 2.1] with a correction from \cite{BFI4}.
\end{remark*}

\section{The dispersion method}\label{sec:dispersion}
In this section we reduce the estimation of $\sum \lambda_q \nu_r \Delta(qr)$ to that of two exponential sums. The reduction is that of [\cite{MaynardLargeModuliI}, Chapter 14], which incorporates the de-amplification by an auxiliary modulus $e$; we need it in a slightly different form.

Let $M,N,Q,R,D,E \geqslant 1$ with $M N \asymp x$, and put
$$
H_1 := \frac{QDR}{M} (\log x)^5, \qquad H_2 := \frac{QDR^2}{M} (\log x)^5.
$$
For complex sequences $(\alpha_n)$, $(\lambda_{d,r})$ and a constant $C \geqslant 0$ we define
\begin{align}
\nonumber \cE_1 := \sum_{e \sim E} \mu^2(e) \sum_{\substack{q \\ (q,a)=1}} \psi_0 \left(\frac{q}{Q} \right) \sum_{\substack{d \sim D \\ (d,a)=1}} & \sum_{\substack{r_1, r_2 \sim R \\ (r_1 r_2,a)=1}} \frac{\lambda_{d, r_1} \ov{\lambda_{d, r_2}}}{\varphi(q d e r_2)\,q d r_1} \\
\nonumber & \times \sum_{\substack{n_1, n_2 \sim N \\ (n_1, q d e r_1)=1 \\ (n_2, q d e r_2)=1}} \alpha_{n_1} \ov{\alpha_{n_2}} \sum_{1 \leqslant |h| \leqslant H_1} \widehat{\psi}_0 \left(\frac{h M}{q d r_1}\right) e\left(\frac{a h \ov{n_1}}{q d r_1}\right)
\end{align}
and
\begin{align}
\nonumber \cE_2:=\sum_{e \sim E} \mu^2(e) \sum_{\substack{q \\ (q,a)=1}} \psi_0 \left(\frac{q}{Q} \right) & \sum_{\substack{d \sim D \\ (d,a)=1}} \sum_{\substack{r_1, r_2 \sim R \\ (r_1,r_2)=1 \\(r_1 r_2,a e)=1}} \frac{\lambda_{d,r_1} \ov{\lambda_{d,r_2}}}{q d r_1 r_2} \sum_{\substack{n_1, n_2 \sim N \\ n_1 \equiv n_2 (\operatorname{mod} qde) \\ (n_1,n_2)=(n_1,r_1)=(n_2,r_2)=1 \\ |n_1-n_2| \geqslant N/(\log x)^{C}}} \alpha_{n_1} \ov{\alpha_{n_2}} \\
\nonumber & \times \sum_{1 \leqslant |h| \leqslant H_2} \widehat{\psi}_0 \left(\frac{h M}{q d r_1 r_2}\right) e\left(a h \frac{n_1-n_2}{q d} \cdot \frac{\ov{n_2 r_1}}{n_1 r_2}\right) e\left(\frac{a h}{q d r_1 r_2 n_1}\right).
\end{align}
Note that in $\cE_2$ the conditions $(n_1,n_2)=1$ and $n_1 \equiv n_2 \pmod{q d e}$ imply $(n_1 n_2,q d e)=1$.

\begin{proposition}[Reduction to exponential sums]\label{prop:dispersion}
Let $\alpha_n, \beta_m, \gamma_{q,d}, \lambda_{d,r}$ be complex sequences with $|\alpha_n|, |\beta_n| \leqslant \tau(n)^{B_0}$, $|\gamma_{q,d}| \leqslant \tau(qd)^{B_0}$ and $|\lambda_{d,r}| \leqslant \tau(dr)^{B_0}$. Suppose that $\alpha_n$ is supported on $P^{-}(n) \geqslant z_0$ and satisfies the Siegel--Walfisz condition \eqref{eq:SW}, and that $\lambda_{d,r}$ is supported on $P^{-}(r) \geqslant z_0$. Let
$$
\cS := \sum_{\substack{d \sim D \\ (d,a)=1}} \sum_{\substack{q \sim Q \\ (q,a)=1}} \sum_{\substack{r \sim R \\ (r,a)=1}} \lambda_{d,r} \gamma_{q,d} \sum_{m \sim M} \beta_m \sum_{n \sim N} \alpha_n \left(\mathbbm{1}_{m n \equiv a (\operatorname{mod} q r d)}-\frac{\mathbbm{1}_{(m n, q r d)=1}}{\varphi(q r d)}\right).
$$
Let $A > 0$, $1 \leqslant E \leqslant x$, and let $C = C(A, B_0)$ be sufficiently large. If
$$
N > QDE(\log x)^{C} \qquad \text{and} \qquad M > (\log x)^{C},
$$
then
$$
|\cS| \ll_{A, B_0} \frac{x}{(\log x)^{A}} + M D^{1/2} Q^{1/2} (\log x)^{O_{B_0}(1)} \left(|\cE_1|^{1/2} + |\cE_2|^{1/2}\right).
$$
\end{proposition}
\begin{proof}
This is [\cite{MaynardLargeModuliI}, Proposition 14.4], in the special case that the coefficients $\lambda_{q,d,r}$ of \cite{MaynardLargeModuliI} do not depend on $q$, except that our sum $\cE_2$ differs from the corresponding sum (denoted $E_2$) of \cite{MaynardLargeModuliI} in three aspects: the condition $(r_1 r_2, e)=1$ is present; the condition $(r_2, q d)=1$ is absent; and the phase is written in a different form. We indicate the necessary changes in the proof.

As in the proof of [\cite{MaynardLargeModuliI}, Proposition 14.4], one inserts the conditions $(n,e)=(r,e)=1$ at a negligible cost, splits the sum over $n$ into residue classes $b \pmod{q d e}$, averages over square-free $e \sim E$ and applies Cauchy's inequality to the sum over $e,b,m,q,d$. This gives
$$
|\cS|^2 \ll M Q D (\log x)^{O(1)} \left(\cS_1 - 2 \operatorname{Re} \cS_2 + \cS_3 \right) + O\left(\frac{x^2}{(\log x)^{2A}}\right),
$$
where $\cS_1, \cS_2, \cS_3$ are the sums evaluated in [\cite{MaynardLargeModuliI}, Lemmas 14.1--14.3] (with $\eta_{q,d,r} = \lambda_{d,r}$). By [\cite{MaynardLargeModuliI}, Lemmas 14.1--14.2], we have
$$
\cS_1 = \cS_{MT} + O\left(\frac{M N^2}{Q D (\log x)^{2A}} \right), \qquad \cS_2 = \cS_{MT} + O(M|\cE_1|) + O\left(\frac{M N^2}{Q D (\log x)^{2A}} \right),
$$
with the main term $\cS_{MT}$ of [\cite{MaynardLargeModuliI}, Lemma 14.1]. It remains to show that
\begin{equation}\label{eq:S3}
\cS_3 = \cS_{MT} + O(M|\cE_2|) + O\left(\frac{M N^2}{Q D (\log x)^{2A}} \right),
\end{equation}
where we recall that
\begin{align}
\nonumber \cS_3 = \sum_{e \sim E} \mu^2(e) \sum_{\substack{q \\ (q,a)=1}} \psi_0 \left(\frac{q}{Q}\right) \sum_{\substack{d \sim D \\ (d,a)=1}} & \sum_{\substack{r_1, r_2 \sim R \\ (r_1 r_2, a e)=1}} \lambda_{d,r_1} \ov{\lambda_{d,r_2}} \\
\nonumber & \times \sum_{\substack{n_1, n_2 \sim N \\ n_1 \equiv n_2 (\operatorname{mod} q d e) \\ (n_1, q d e r_1) = (n_2, q d e r_2)=1}} \alpha_{n_1} \ov{\alpha_{n_2}} \sum_{\substack{m \\ m n_1 \equiv a (\operatorname{mod} q d r_1) \\ m n_2 \equiv a (\operatorname{mod} q d r_2)}} \psi_0 \left(\frac{m}{M}\right).
\end{align}
We follow the proof of [\cite{MaynardLargeModuliI}, Lemma 14.3]. Exactly as there, one may restrict the summation to $(r_1, r_2)=1$, to $(n_1, n_2)=1$ and to $|n_1 - n_2| \geqslant N/(\log x)^C$, at the cost of an error $O(M N^2/(Q D (\log x)^{2A}))$. We deviate from \cite{MaynardLargeModuliI} in two places:

\noindent\textit{(1)} We do not remove the condition $(e, r_1 r_2)=1$. In \cite{MaynardLargeModuliI} this condition is removed at the end of the first part of the proof and re-introduced later in the main term; neither step is needed.

\noindent\textit{(2)} We do not discard the terms with $(q d, r_2)>1$. In \cite{MaynardLargeModuliI} this is done only in order to apply the Chinese remainder theorem in the form of Lemma~\ref{lem:bezout}; we use Lemma~\ref{lem:reciprocity} instead. Indeed, after the above restrictions the hypotheses of Lemma~\ref{lem:reciprocity} are satisfied with $\ell = q d$ and $b = a$, so the two congruences for $m$ are equivalent to $m \equiv \mu \pmod{q d r_1 r_2}$, and Poisson summation [\cite{MaynardLargeModuliI}, Lemma 13.4] gives
$$
\sum_{m \equiv \mu (\operatorname{mod} q d r_1 r_2)} \psi_0 \left(\frac{m}{M}\right) = \frac{M \widehat{\psi}_0(0)}{q d r_1 r_2} + \frac{M}{q d r_1 r_2} \sum_{1 \leqslant |h| \leqslant H_2} \widehat{\psi}_0 \left(\frac{h M}{q d r_1 r_2} \right) e\left(\frac{h \mu}{q d r_1 r_2} \right) + O(x^{-100}),
$$
where, by Lemma~\ref{lem:reciprocity},
$$
e\left(\frac{h \mu}{q d r_1 r_2} \right) = e\left(a h \frac{n_1 - n_2}{q d} \cdot \frac{\ov{n_2r_1}}{n_1r_2} \right) e\left(\frac{a h}{q d r_1 r_2 n_1}\right).
$$
The terms with $h \neq 0$ contribute $M \cE_2$ to $\cS_3$. In the contribution of the terms with $h = 0$ one inserts, as in \cite{MaynardLargeModuliI}, the factor
$$
\frac{\varphi(q d e)^2 \varphi(q d r_1 r_2)}{\varphi(q d) \varphi(q d e r_1) \varphi(q d e r_2)} = \frac{r_1 r_2 \prod_{p \mid r_1 r_2,\ p \nmid q d}\frac{p - 1}{p}}{r_1 \prod_{p \mid r_1,\ p \nmid q d e}\frac{p - 1}{p} \cdot r_2 \prod_{p \mid r_2,\ p \nmid q d e}\frac{p - 1}{p}},
$$
which in \cite{MaynardLargeModuliI} is $1+O(\log x/z_0)$ and in our situation is exactly equal to $1$, since $(r_1, r_2) = (r_1 r_2, e)=1$. From this point on the argument of \cite{MaynardLargeModuliI} applies verbatim: the conditions $(r_1, r_2)=1$, $(n_1, n_2)=1$ and $|n_1 - n_2| \geqslant N/(\log x)^C$ are removed from the terms with $h = 0$ at the cost of an admissible error, and the condition $n_1 \equiv n_2 \pmod{q d e}$ is replaced by the factor $\frac{1}{\varphi(q d e)}$ by means of the Siegel--Walfisz condition and the large sieve [\cite{MaynardLargeModuliI}, Lemma 13.8]; here the hypothesis $N > Q D E (\log x)^C$ is used. This gives \eqref{eq:S3}, and Proposition~\ref{prop:dispersion} follows.
\end{proof}

\section{The first exponential sum}\label{sec:E1}
The sum $\cE_1$ is easily handled, because the sum over $q$ is an incomplete Kloosterman sum with a smooth weight. The following lemma is a variant of [\cite{MaynardLargeModuliII}, Lemma 5.10].
\begin{lemma}\label{lem:E1}
Let $\eta > 0$ and $1 \leqslant E \leqslant x$, and suppose that
\begin{equation}\label{eq:E1cond}
D R N^{3/2} \leqslant x^{1 - 2 \eta} \qquad \text{and} \qquad Q D R \leqslant x^{1 - 2 \eta}.
\end{equation}
Then for all $1$-bounded sequences $(\alpha_n)$ and $(\lambda_{d,r})$, we have
$$
\cE_1 \ll_\eta \frac{N^2}{Q D x^{\eta}}.
$$
\end{lemma}
\begin{proof}
Fix $e,d,r_1,r_2,n_1,n_2,h$ with $(n_1, d e r_1)=(n_2, d e r_2)=1$ and put $m_0 := d e r_2$. The sum over $q$ in $\cE_1$ is
$$
\Sigma := \sum_{(q,a n_1 n_2)=1} \psi_0 \left(\frac{q}{Q}\right) \frac{1}{q \varphi(q m_0)} \widehat{\psi}_0 \left(\frac{h M}{q d r_1} \right) e\left(\frac{a h \ov{n_1}}{q d r_1}\right).
$$
By Lemma~\ref{lem:bezout} with $(u,v)=(n_1,q d r_1)$ we have
$$
e\left(\frac{a h \ov{n_1}}{q d r_1}\right) = e\left(-\frac{a h \ov{q d r_1}}{n_1} \right) e\left(\frac{a h}{q d r_1 n_1}\right).
$$
Moreover
$$
\frac{1}{\varphi(q m_0)} = \frac{1}{q \varphi(m_0)} \sum_{\substack{g \mid q \\ (g, m_0)=1}} \frac{\mu^2(g)}{\varphi(g)}, \qquad \mathbbm{1}_{(q,a n_2)=1} = \sum_{t \mid (q, a n_2)} \mu(t).
$$
Hence, writing $q = [g,t] q'$,
$$
\Sigma = \frac{1}{\varphi(m_0)} \sum_{(g, m_0 n_1)=1} \frac{\mu^2(g)}{\varphi(g)} \sum_{\substack{t \mid a n_2 \\ (t, n_1)=1}} \mu(t) \sum_{(q', n_1)=1} W\left([g,t] q'\right) e\left(\frac{b \ov{q'}}{n_1}\right),
$$
where $b := -a h \ov{d r_1 [g,t]} \pmod{n_1}$ and
$$
W(q) := \psi_0 \left(\frac{q}{Q}\right) \frac{1}{q^2} \widehat{\psi}_0 \left(\frac{h M}{q d r_1}\right) e\left(\frac{a h}{q d r_1 n_1}\right).
$$
For $q \asymp Q$ we have $|h M/(q d r_1)| \ll H_1 M/(QDR) = (\log x)^5$ and $|a h /(q d r_1 n_1)| \ll 1$, so that the function $f(u) := Q^2 W(Q u)$ is supported on $[1/2, 5/2]$ and satisfies $\|f^{(j)}\|_\infty \ll_j (\log x)^{5j}$. Since $(b, n_1)=(ah, n_1)$, Lemma~\ref{lem:weil} (with $Y = Q/[g,t]$) shows that the innermost sum is
$$
\ll x^{o(1)} Q^{-2} \left(N^{1/2} + \frac{Q}{N}(a h,n_1)\right).
$$
Summing over $g \leqslant 3 Q$ and over the divisors $t$ of $a n_2$ we obtain 
$$
\Sigma \ll x^{o(1)} \frac{1}{\varphi(m_0) Q^{2}} \left(N^{1/2} + Q N^{-1}(a h, n_1)\right).
$$
Therefore, using $\sum_{1 \leqslant |h| \leqslant H_1}(a h, n_1) \ll H_1 \tau(n_1)$,
\begin{align}
\nonumber \cE_1 \ll&\ x^{o(1)} \sum_{e \sim E} \sum_{d \sim D} \sum_{r_1, r_2 \sim R} \frac{1}{d r_1 \varphi(d e r_2)} \sum_{n_1, n_2 \sim N} \sum_{1 \leqslant |h| \leqslant H_1} \frac{1}{Q^{2}} \left(N^{1/2} + \frac{Q}{N}(a h,n_1)\right) \\
\nonumber \ll&\ x^{o(1)} \frac{N^2 H_1}{D Q^2} \left(N^{1/2} + \frac{Q}{N}\right) \\
\nonumber =&\ x^{o(1)} \frac{N^2}{Q D} \cdot \frac{D R N^{3/2} + Q D R}{M N},
\end{align}
and the lemma follows from \eqref{eq:E1cond} and $M N \asymp x$.
\end{proof}

\section{The second exponential sum}\label{sec:E2}
In this section we prove an upper bound for $\cE_2$, which is the most crucial estimate in the proof of Theorem~\ref{thm:typeII} and Theorem~\ref{thm:main}.
\begin{proposition}\label{prop:E2}
Let $\vartheta$ be admissible. Let $\eta > 0$ and $C \geqslant 0$ be fixed. Let $M,N,Q,R,D,E \geqslant 1$ satisfy $M N \asymp x$ and
\begin{align}
& Q \leqslant M, \qquad D \leqslant y_0, \qquad Q D E \leqslant N, \label{eq:P0} \\
& N D \leqslant x^{-3\eta} M E, \label{eq:P1} \\
& E^{7/2} \left(1 + \frac{N}{M}\right)^{1/2} N^{1/2} Q^{1/2} R^{4} \leqslant x^{1-4\eta}, \label{eq:P2} \\
& E^{3} \left(1 + \frac{N}{M}\right)^{1/2} N^{1/2 + \vartheta} Q^{1/2} R^{4 - 3 \vartheta} \leqslant x^{1-4\eta}. \label{eq:P3}
\end{align}
Let $(\alpha_n)$ and $(\lambda_{d,r})$ be $1$-bounded complex sequences such that $\lambda_{d,r}$ is supported on $P^{-}(r) \geqslant z_0$. Then
$$
\cE_2 \ll_{a, \eta} \frac{N^2}{Q D x^{\eta}}.
$$
\end{proposition}
The rest of this section is devoted to the proof. We note that \eqref{eq:P2} implies $R \leqslant x^{1/4}$, and that $H_2 = Q D R^2 (\log x)^5/M \leqslant D R^2 (\log x)^5 \leqslant x^{1/2+o(1)}$ by \eqref{eq:P0}.

\subsection{Passage to the complementary divisor}
The trivial bound for $\cE_2$ is
$$
\sum_{e \sim E} \sum_{d \sim D} \sum_{Q/2 \leqslant q \leqslant 5Q/2} \frac{1}{q d} \sum_{r_1, r_2 \sim R} \frac{1}{r_1 r_2} \sum_{\substack{n_1, n_2 \sim N \\ n_1 \equiv n_2 (\operatorname{mod} q d e)}} \sum_{1 \leqslant |h| \leqslant H_2} 1 \ll \frac{N^2 H_2}{Q D},
$$
since $Q D E \leqslant N$. In $\cE_2$ we have 
$$
\left|\frac{a h}{q d r_1 r_2 n_1}\right| \ll \frac{H_2}{Q D R^2 N} \ll \frac{(\log x)^5}{x}.
$$
Hence the factor $e(a h/(q d r_1 r_2 n_1))$ may be replaced by $1$ at the cost of an error
$$
O\left(\frac{N^2 H_2 (\log)^5}{Q D x}\right) = O\left(\frac{N^2}{Q D x^{1/3}}\right).
$$

Next, for $d, e, n_1, n_2$, as in $\cE_2$ we have $n_1 \neq n_2$, and the variable $q$ runs over positive integers such that $q d e$ divides $n_1 - n_2$. We write
$$
n_1 - n_2 = q,\ d,\ e,\ f, \qquad f \in \Z \backslash \{0\},
$$
and use $f$ as a variable in place of $q$; thus $q = (n_1 - n_2)/(d e f)$. Since $q \geqslant Q/2$, $d > D$, $e > E$ and $|n_1 - n_2| < N$, we have
$$
1 \leqslant |f| \leqslant F := \frac{2N}{Q D E},
$$
and $(n_1 - n_2)/(q d) = e f$. Therefore
\begin{equation}\label{eq:E2prime}
\cE_2 = \cE_2' + O\left(\frac{N^2}{Q D x^{1/3}}\right),
\end{equation}
where
\begin{align}
\nonumber \cE_2' := \sum_{e \sim E} & \mu^2(e) \sum_{\substack{d \sim D \\ (d,a)=1}} \frac{1}{d} \sum_{1 \leqslant |f| \leqslant F} \sum^{\mathbf{N}}_{n_1, n_2} \alpha_{n_1} \ov{\alpha_{n_2}} \\
\nonumber & \times \sum^{\mathbf{R}}_{r_1,r_2} \frac{\lambda_{d, r_1} \ov{\lambda_{d,r_2}}}{r_1 r_2} \sum_{1 \leqslant |h| \leqslant H_2} \frac{1}{q} \psi_0 \left(\frac {q}{Q}\right) \widehat{\psi}_0 \left(\frac{h M}{q d r_1 r_2}\right) e\left(a h e f \frac{\ov{n_2r_1}}{n_1r_2}\right), \qquad q := \frac{n_1 - n_2}{d e f}.
\end{align}
Here $\sum^\mathbf{N}$ denotes a sum over $n_1, n_2 \sim N$ with
\begin{equation}\label{eq:sum_n}
(n_1, n_2) = 1,\ n_1 \equiv n_2 (\operatorname{mod} d e f),\ \frac{n_1 - n_2}{d e f} > 0,\ \left(\frac{n_1 - n_2}{d e f}, a \right) = 1,\ |n_1 - n_2| \geqslant \frac{N}{(\log x)^C},
\end{equation}
and $\sum^\mathbf{R}$ denotes a sum over $r_1, r_2 \sim R$ with
\begin{equation}\label{eq:sum_r}
(r_1, r_2) = 1, \qquad (r_1 r_2, a e) = 1, \qquad (r_1, n_1) = (r_2, n_2) = 1.
\end{equation}
Note that \eqref{eq:sum_n} implies $(n_1 n_2, e) = 1$.

\subsection{Separation of variables and Cauchy's inequality}
The variable $q$, which is now a function of $n_1,n_2,d,e,f$, is entangled with $h,r_1,r_2$ through the factor $\widehat{\psi}_0$. We separate them as follows. By the definition of the Fourier transform and the change of variables $t = q \tau$,
$$
\widehat{\psi}_0 \left(\frac{h M}{q d r_1 r_2}\right) = \int \psi_0(t) e\left(-\frac{t h M}{q d r_1 r_2}\right) \dd t = q \int \psi_0(q \tau) e\left(-\frac{\tau h M}{d r_1 r_2}\right) \dd \tau.
$$
Since $\psi_0(q/Q) \psi_0(q \tau)$ vanishes unless $1/(5Q) \leqslant \tau \leqslant 5/Q$, we get
\begin{equation}\label{eq:separation}
\frac{1}{q} \psi_0\left(\frac{q}{Q}\right) \widehat{\psi}_0 \left(\frac{h M}{q d r_1 r_2}\right) = \int_{1/(5Q)}^{5/Q} \psi_0 \left(\frac{q}{Q}\right) \psi_0(q \tau) e\left(-\frac{\tau h M}{d r_1 r_2} \right) \dd \tau.
\end{equation}
For $\tau \in [1/(5Q), 5/Q]$ and $d \sim D$ we define the $1$-bounded coefficients
$$
\xi(h, r_1, r_2) = \xi_{\tau,d}(h, r_1, r_2) := \frac{R^2 \lambda_{d,r_1} \ov{\lambda_{d,r_2}}}{r_1 r_2} e\left(-\frac{\tau h M}{d r_1 r_2}\right),
$$
which are supported on $r_1, r_2 \sim R$ with $P^{-}(r_1r_2) \geqslant z_0$. Let $\cN$ be the set of pairs $(r_1,r_2)$ with $r_1, r_2 \sim R$, $(r_1, r_2) = 1$ and $P^{-}(r_1 r_2) \geqslant z_0$, and let $\cN = \cN_1 \sqcup \cdots \sqcup \cN_J$ be the partition given by Lemma~\ref{lem:coprime}. For $1 \leqslant j \leqslant J$ and $e,f,n_1,n_2$ with $(n_1, n_2) = 1$ put
$$
I_j(f, n_1, n_2) := \sum_{\substack{(r_1, r_2) \in \cN_j \\ (r_1 r_2, a e) = 1 \\ (r_1, n_1) = (r_2, n_2) = 1}} \sum_{1 \leqslant |h| \leqslant H_2} \xi(h, r_1, r_2) e\left(a h e f \frac{\ov{n_2 r_1}}{n_1 r_2}\right).
$$
Inserting \eqref{eq:separation} into the definition of $\cE_2'$ we obtain
$$
|\cE_2'| \leqslant \frac{1}{R^2} \int_{1/(5Q)}^{5/Q} \sum_{e \sim E} \mu^2(e) \sum_{d \sim D} \frac{1}{d} \sum_{j \leqslant J} \sum_{1 \leqslant |f| \leqslant F} \sum^{\mathbf{N}}_{n_1, n_2} \|\psi_0\|_\infty^2 \left|I_j(f, n_1, n_2)\right| \dd \tau.
$$
For fixed $d$ and $e$ the number of triples $(f, n_1, n_2)$ with $1 \leqslant |f| \leqslant F$ and $n_1, n_2$ satisfying \eqref{eq:sum_n} with $Q/2 \leqslant q \leqslant 5Q/2$ is $\ll N^2/(D E)$, because such a triple is determined by $n_1 \sim N$, $q \in [Q/2, 5Q/2]$ and $f$. By Cauchy's inequality, and enlarging the range of summation by positivity, we infer that
\begin{equation}\label{eq:afterCS}
\cE_2' \ll \frac{J}{R^2} \cdot \frac{1}{Q} \cdot E \cdot D \cdot \frac{1}{D} \cdot \left(\frac{N^2}{D E}\right)^{1/2} \sup \cZ^{1/2} = \frac{J N E^{1/2}}{Q D^{1/2} R^2} \sup \cZ^{1/2},
\end{equation}
where the supremum is over $\tau \in [1/(5Q), 5/Q]$, $d \sim D$, square-free $e \sim E$ and $j \leqslant J$, and
\begin{equation}\label{eq:Zdef}
\cZ = \cZ(\tau,d,e,j) := \sum_{1 \leqslant |f| \leqslant F} \sum_{\substack{n_1, n_2 \\ (n_1, n_2) = 1 \\ n_1 \equiv n_2 (\operatorname{mod} e) \\ (n_1 n_2, e) = 1}} \Psi\left(\frac{n_1}{N}\right) \Psi\left(\frac{n_2}{N}\right) \left|I_j(f, n_1, n_2)\right|^2.
\end{equation}
In passing from \eqref{eq:sum_n} to the conditions in \eqref{eq:Zdef} we have dropped the congruence $n_1 \equiv n_2 \pmod{d f}$ and the conditions involving $q$, but we have kept the congruence modulo $e$. We shall prove that
\begin{equation}\label{eq:Ztarget}
\cZ \ll \frac{R^4 N^2}{D E x^{5\eta/2 -o(1)}}
\end{equation}
uniformly in $\tau,d,e,j$. By \eqref{eq:afterCS}, \eqref{eq:E2prime} and $J = x^{o(1)}$ this gives $\cE_2 \ll N^2/(QD x^{5\eta/4 -o(1)})$, which proves Proposition~\ref{prop:E2}.

\subsection{The diagonal terms}
We expand the square in \eqref{eq:Zdef}. Let $(h, r_1, r_2)$ and $(h', r_1', r_2')$ be two triples occurring in $I_j(f, n_1, n_2)$. Since $(r_1, r_2)$ and $(r_1', r_2')$ belong to the same class $\cN_j$ we have $(r_1, r_2') = (r_1', r_2) = 1$ in addition to $(r_1, r_2) = (r_1', r_2') = 1$; hence
$$
(r_1 r_1', r_2 r_2') = 1.
$$
Together with $(n_1, n_2) = (n_1, r_1 r_1') = (n_2, r_2 r_2') = 1$ this gives $(n_2 r_1 r_1', n_1 r_2 r_2') = 1$, and Lemma~\ref{lem:combine} yields
\begin{align}\label{eq:Zexpanded}
\nonumber \cZ = \sum_{\substack{(h, r_1, r_2) \\ (h', r_1', r_2')}} & \xi(h, r_1, r_2) \ov{\xi(h', r_1', r_2')} \\
& \times \sum_{1 \leqslant |f| \leqslant F} \sum_{\substack{n_1, n_2 \\ (n_1, n_2) = 1 \\ n_1 \equiv n_2 (\operatorname{mod} e), (n_1 n_2, e) = 1 \\ (n_1, r_1 r_1') = (n_2, r_2 r_2') = 1}} \Psi\left(\frac{n_1}{N}\right) \Psi\left(\frac{n_2}{N}\right) e\left(a e f k \frac{\ov{n_2 \rho}}{n_1 \sigma}\right),
\end{align}
where
$$
\rho := r_1 r_1', \qquad \sigma := [r_2, r_2'], \qquad k := \frac{h r_1' r_2' - h' r_1 r_2}{(r_2, r_2')},
$$
and the outer sum is over all $h, h'$ with $1 \leqslant |h|, |h'| \leqslant H_2$ and over all pairs $(r_1, r_2)$, $(r_1', r_2')$ in $\cN_j$ which satisfy $(r_1 r_2 r_1' r_2', a e) = 1$. We write $\cZ = \cZ_{=} + \cZ_{\neq}$ according to whether $k = 0$ or $k \neq 0$. If $k = 0$, then $h r_1' r_2' = h' r_1 r_2$. Given $(h, r_1', r_2')$, the integer $h' r_1 r_2$ is determined, and so $(h', r_1, r_2)$ is determined up to $x^{o(1)}$ possibilities by the divisor bound. Thus there are at most $x^{o(1)} H_2 R^2$ pairs of triples with $k = 0$. For each of them the sum over $f, n_1, n_2$ in \eqref{eq:Zexpanded} is trivially $\ll F N (N/E + 1) \ll F N^2/E$. Hence
$$
\cZ_{=} \ll x^{o(1)} H_2 R^2 \frac{F N^2}{E} = x^{o(1)} \frac{R^4 N^3}{M E^2},
$$
because $H_2 F = 2 R^2 N (\log)^5/(M E)$. By \eqref{eq:P1} we have $N/(M E)\le x^{-3\eta}/D$, so that
\begin{equation}\label{eq:Zdiag}
\cZ_{=} \ll \frac{R^4 N^2}{D E x^{3\eta -o(1)}}.
\end{equation}
This is the only place where the auxiliary modulus $e$ is essential: for $E = 1$ the bound \eqref{eq:Zdiag} would require $N \leqslant x^{-3\eta} M$.

\subsection{The off-diagonal terms}
We now estimate $\cZ_{\neq}$ by means of Theorem~\ref{thm:DI}. We split the summation in \eqref{eq:Zexpanded} according to the sign $\vartheta \in \{\pm{1}\}$ of $a f k$ and according to the sizes of $\rho$ and $\sigma$. Since $R^2 < \rho \leqslant 4R^2$ and $R < \sigma \leqslant 4R^2$, it suffices to consider $O(\log x)$ pairs $(R^{*}, S^{*})$ with
$$
R^2 \leqslant R^{*} \leqslant 2R^2, \qquad R \leqslant S^{*} < 4R^2, \qquad R^{*} < \rho \leqslant 2R^{*}, \qquad S^{*} < \sigma \leqslant 2S^{*}.
$$
Fix $\vartheta$, $R^{*}$ and $S^{*}$. In the corresponding terms we have $(r_2, r_2') = r_2 r_2'/\sigma > R^2/(2S^{*})$, so that $|k| \leqslant 8 H_2 R^2/(r_2, r_2') < 16 H_2 S^{*}$, and therefore
$$
1 \leqslant |a e f k| \leqslant N^{*} := 32|a| E F H_2 S^{*} = 64|a| (\log x)^5 \frac{N}{M} R^2 S^{*}.
$$
For integers $n \geqslant 1$, $\rho \sim R^{*}$, $\sigma \sim S^{*}$ we define
$$
b_{n,\rho,\sigma} := \sum\xi(h,r_1,r_2) \ov{\xi(h',r_1',r_2')},
$$
where the sum is over all $f$ with $1 \leqslant |f| \leqslant F$ and all pairs of triples as in \eqref{eq:Zexpanded} such that
$$
r_1 r_1' = \rho, \qquad [r_2, r_2'] = \sigma, \qquad a e f k = \vartheta n.
$$
On the support of $b_{n,\rho,\sigma}$ we have $(e, \sigma) = (\rho, \sigma) = 1$. Hence, for $(n_1, e) = 1$, the conditions $(n_1, n_2) = (n_1, r_1 r_1') = (n_2, r_2 r_2') = 1$ in \eqref{eq:Zexpanded} are equivalent to the single condition $(e \rho n_2, \sigma n_1) = 1$. Splitting $n_1 \equiv n_2 \pmod{e}$ into the reduced residue classes $u$ modulo $e$, we see that the part of $\cZ_{\neq}$ under consideration equals
\begin{equation}\label{eq:DIshape}
\sum_{\substack{u (\operatorname{mod} e) \\ (u, e) = 1}} \sum_{n \leqslant N^{*}} \sum_{\rho \sim R^{*}} \sum_{\sigma \sim S^{*}} b_{n,\rho,\sigma} \sum_{\substack{n_1 \equiv u, n_2\equiv u (\operatorname{mod} e) \\ (e \rho n_2, \sigma n_1) = 1}} \Psi\left(\frac{n_1}{N}\right) \Psi\left(\frac{n_2}{N}\right) e\left(\vartheta n \frac{\ov{\rho n_2}}{\sigma n_1}\right).
\end{equation}
If $\vartheta = -1$ we pass to the complex conjugate, which replaces $b$ by $\ov{b}$. Let $\Psi = \Psi_1 + \Psi_2 + \Psi_3$ be a smooth partition with $\Psi_i$ supported in $[c_i, 1.9 c_i]$ for $c_1 = 0.45$, $c_2 = 0.8$, $c_3 = 1.4$, and let $w(n,\rho,\sigma)$ be a smooth function, compactly supported in $\R_{>0}^3$, which is equal to $1$ on $[1, N^{*}] \times [R^{*}, 2R^{*}] \times [S^{*}, 2S^{*}]$ and satisfies $\partial^{\nu}w \ll_\nu n^{-\nu_3} \rho^{-\nu_4} \sigma^{-\nu_5}$. Then, for each $u$ and each pair $(i, i')$, the inner part of \eqref{eq:DIshape} is a sum of the type considered in Theorem~\ref{thm:DI}, with
$$
(c,d,n,r,s) = (n_1,n_2,n,\rho,\sigma), \qquad g = \Psi_i\left(\frac {c}{N}\right) \Psi_{i'}\left(\frac{d}{N}\right) w(n,r,s),
$$
$$
(\sC,\sD,\sN,\sR,\sS) = (0.99 c_iN, 0.99 c_{i'}N, N^{*}, R^{*}, S^{*}), \qquad \sq = e, \qquad c_0 = d_0 = u.
$$
We now compute the quantity $\sK$ of \eqref{eq:K}. Using $R^{*} \leqslant 2R^2$, $S^{*} < 4R^2$ and $e \leqslant 2E$, we find
$$
\sN + \sR \sS \ll (\log x)^5 \left(1 + \frac{N}{M}\right) R^2 S^{*}, \qquad \sC + \sD \sR \ll N R^2,
$$
and hence
\begin{align}
\nonumber \sq (\sC \sS (\sN + \sR \sS) (\sC + \sD \sR) + \sN \sR \sS) \ll&\ (\log x)^5 \left(1 + \frac{N}{M}\right) E N^2 R^8, \\
\nonumber \sD^2 \sN \sR \ll&\ (\log x)^5 \left(1 + \frac{N}{M}\right) N^2 R^6, \\
\nonumber \sC \sS (\sC \sD \sqrt{\sR})^{2 \vartheta} (\sN + \sR \sS)^{1 - \vartheta} (\sC + \sD \sR)^{1 - 2 \vartheta} \ll&\ (\log x)^5 \left(1 + \frac{N}{M}\right) N^{2 + 2 \vartheta} R^{8 - 6 \vartheta}.
\end{align}
The second term is dominated by the first, and so
\begin{equation}\label{eq:Kbound}
\sK \ll (\log x)^{5/2} \left(1 + \frac{N}{M}\right)^{1/2} \left(E^{1/2} N R^4 + N^{1 + \vartheta} R^{4 - 3 \vartheta}\right).
\end{equation}

Next we estimate $\|b\|_2$. Write $\mathbf{r} := (r_1, r_1', r_2, r_2')$. Since $\xi$ is $1$-bounded, $|b_{n, \rho, \sigma}|$ is at most the number $\nu(n, \rho, \sigma)$ of tuples $(f, h, h', \mathbf{r})$ counted in its definition. Given $(n, \rho, \sigma)$ there are at most $x^{o(1)}$ possibilities for $(f, k, \mathbf{r})$, because $f$ divides $n$, $r_1 r_1' = \rho$, and $r_2, r_2'$ divide $\sigma$. For given $(k, \mathbf{r})$ let $\nu_0(k; \mathbf{r})$ be the number of pairs $(h, h')$ with $1 \leqslant |h|,|h'| \leqslant H_2$ and
$$
h r_1' r_2' - h' r_1 r_2 = k (r_2, r_2').
$$
The solutions $(h, h') \in \Z^2$ of this linear equation form a coset of the lattice generated by the vector $(r_1 r_2/g, r_1' r_2'/g)$, where
$$
g := (r_1 r_2, r_1'r_2') = (r_1, r_1') (r_2, r_2')
$$
in view of $(r_1 r_1', r_2 r_2') = 1$. As $r_1 r_2 > R^2$, we get
$$
\nu_0(k; \mathbf{r}) \leqslant 1 + \frac{2 H_2(r_1, r_1')(r_2, r_2')}{R^2}.
$$
By Cauchy's inequality $\nu(n, \rho, \sigma)^2 \leqslant x^{o(1)} \sum \nu_0(k; \mathbf{r})^2$, the sum being over the tuples $(f, k, \mathbf{r})$ compatible with $(n, \rho, \sigma)$. Summing over $(n, \rho, \sigma)$ and using that $\sum_k \nu_0(k; \mathbf{r}) \leqslant (2 H_2)^2$, we obtain
\begin{align}
\nonumber \|b\|_2^2 \ll&\ x^{o(1)} \sum_{1 \leqslant |f| \leqslant F} \sum_{\substack{r_1, r_1', r_2, r_2' \sim R \\ (r_1 r_1', r_2 r_2') = 1}} \left(\max_k \nu_0(k; \mathbf{r})\right) \sum_k \nu_0(k; \mathbf{r}) \\
\nonumber \ll&\ x^{o(1)} F H_2^2 \sum_{\substack{r_1, r_1', r_2, r_2' \sim R \\ (r_1 r_1', r_2 r_2') = 1}} \left(1 + \frac{H_2(r_1, r_1')(r_2, r_2')}{R^2}\right) \\
\nonumber \leqslant&\ x^{o(1)} F H_2^2 \sum_{r_1, r_1', r_2, r_2' \sim R} \left(1 + \frac{H_2(r_1, r_1')(r_2, r_2')}{R^2}\right) \\
\nonumber \ll&\ x^{o(1)} F H_2^2 R^4 \left(1 + \frac{H_2}{R^2}\right),
\end{align}
since $\sum_{r, r' \sim R}(r, r') \ll R^2(\log x)$. By \eqref{eq:P0} we have $H_2/R^2 = Q D (\log x)^5/M \leqslant x^{o(1)}$. Inserting the values of $F$ and $H_2$ we conclude that
\begin{equation}\label{eq:bnorm}
\|b\|_2 \ll x^{o(1)} F^{1/2} H_2 R^2 = x^{o(1)} \frac{(N Q D)^{1/2} R^4}{M E^{1/2}}.
\end{equation}

We apply Theorem~\ref{thm:DI} with $\epsilon$ and $\epsilon_0$ sufficiently small such that $(\sq\sC\sD\sN\sR\sS)^{\epsilon +O(\epsilon_0)} \leqslant x^{\eta/2}$. Summing over the $\varphi(e) \leqslant 2E$ classes $u$, over $(i, i')$, over $\vartheta$ and over the $O(\log x)$ pairs $(R^{*},S^{*})$, we deduce from \eqref{eq:Kbound} and \eqref{eq:bnorm} that
\begin{align}
\nonumber \cZ_{\neq} \ll&\ x^{\eta/2 +o(1)} E^{5/2} \left(1+\frac{N}{M}\right)^{1/2} \left(E^{1/2} N R^4 + N^{1 + \vartheta} R^{4 - 3 \vartheta}\right) \frac{(N Q D)^{1/2} R^4}{M E^{1/2}} \\
\nonumber =&\ x^{\eta/2+o(1)} \frac{R^4 N^2}{D E} \cdot \frac{D^{3/2}(1+N/M)^{1/2}}{M N} \left(E^{7/2} N^{1/2} Q^{1/2} R^4 + E^{3} N^{1/2 + \vartheta} Q^{1/2} R^{4 - 3 \vartheta}\right).
\end{align}
By \eqref{eq:P2}, \eqref{eq:P3}, $D \leqslant y_0$ and $M N \asymp x$, the last factor is $\ll x^{-4\eta +o(1)}$. Hence
$$
\cZ_{\neq} \ll \frac{R^4 N^2}{D E x^{7\eta/2 -o(1)}},
$$
and together with \eqref{eq:Zdiag} this proves \eqref{eq:Ztarget}. The proof of Proposition~\ref{prop:E2} is now complete.

\section{Proof of Theorem \ref{thm:typeII} and Corollary \ref{cor:symmetric}}\label{sec:proofTypeII}
\begin{proof}[Proof of Theorem~\ref{thm:typeII}]
Let $\cS_0$ denote the sum to be estimated. Note that \eqref{eq:T1} and $M N \asymp x$ give
\begin{equation}\label{eq:NM}
\frac{N}{M} \ll x^{\varepsilon/50}, \qquad Q \leqslant N x^{-\varepsilon} \leqslant M, \qquad R \leqslant x^{1/4}, \qquad Q R \leqslant x^{3/4 +\varepsilon/100},
\end{equation}
where the bound for $R$ follows from \eqref{eq:T2}. If $R \leqslant x^{\varepsilon/2}$, then all moduli $q r$ are at most $4 Q R \leqslant 4 N x^{-\varepsilon/2} \leqslant x^{1/2 -\varepsilon/3}$, and the claim follows from the generalized Bombieri--Vinogradov theorem (see [\cite{BFI}, Theorem 0]). We may therefore assume that $R > x^{\varepsilon/2}$. We write each $r \sim R$ uniquely as $r = d r'$ with $d = \sm(r; z_0)$, so that $P^{+}(d) \leqslant z_0 < P^{-}(r')$. The moduli $q r$ with $d \geqslant y_0$ satisfy $\sm(q r; z_0) \geqslant y_0$, and since the coefficient $\sum_{q r = n} \lambda_q \nu_r$ of a modulus $n$ is bounded by $\tau(n)^{2B_0 +1}$, [\cite{MaynardLargeModuliI}, Lemma 13.10] shows that these moduli contribute $O(x/(\log x)^A)$ to $\cS_0$. We split the remaining range of $d$ into $O(\log x)$ dyadic intervals $d \sim D$, $D \leqslant y_0$ (the same results also hold for $d = 1$); then $r' \sim R'$ for some $R/(2D) \leqslant R' \leqslant R/D$ and $R' \geqslant 1$. For each such $(D, R')$, the corresponding part of $\cS_0$ is a sum $\cS$ as in Proposition~\ref{prop:dispersion}, with $R'$ in place of $R$ and with
$$
\gamma_{q,d} := \lambda_q, \qquad \lambda_{d,r'} := \nu_{d r'} \mathbbm{1}_{P^{+}(d) \leqslant z_0<P^{-}(r')} \mathbbm{1}_{R < d r' \leqslant 2R}.
$$
It suffices to show that each of these sums is $O(x/(\log x)^{A+1})$.

We apply Proposition~\ref{prop:dispersion} with $E := x^{\varepsilon/10}$. Its hypotheses hold, since $N \geqslant Q x^{\varepsilon} > Q D E (\log x)^{C}$ and $M \gg x^{1/2 -\varepsilon/100}$. As $|\alpha_n|, |\lambda_{d,r'}| \leqslant \tau(\cdot)^{B_0} \leqslant x^{o(1)}$, it remains to verify that Lemma~\ref{lem:E1} and Proposition~\ref{prop:E2} (applied to suitably normalized sequences, with $R'$ in place of $R$ and with $\eta := \varepsilon/100$) yield
\begin{equation}\label{eq:E1E2}
\cE_1, \cE_2 \ll \frac{N^2}{Q D x^{\eta -o(1)}},
\end{equation}
because then Proposition~\ref{prop:dispersion} gives $\cS \ll x/(\log x)^{A+1} + M N x^{-\eta/2 +o(1)}$.

For Lemma~\ref{lem:E1} we note that by \eqref{eq:T1}, \eqref{eq:T2} and \eqref{eq:NM}
\begin{align}
\nonumber D R' N^{3/2} \leqslant&\ R N^{3/2} \leqslant \left(x^{1-\varepsilon} N^{-1/2} \right)^{1/4} N^{3/2} \leqslant x^{1/4} N^{11/8} \leqslant x^{15/16 +\varepsilon/50}, \\
\nonumber Q D R' \leqslant&\ Q R \leqslant x^{3/4 +\varepsilon/100},
\end{align}
so that \eqref{eq:E1cond} holds. For Proposition~\ref{prop:E2} we check \eqref{eq:P0}--\eqref{eq:P3}. Condition \eqref{eq:P0} follows from \eqref{eq:NM} and $Q D E \leqslant Q y_0 x^{\varepsilon/10} \leqslant N$. By \eqref{eq:NM},
$$
N D \leqslant x^{\varepsilon/50 +o(1)} M \leqslant x^{-3\eta} M E,
$$
which is \eqref{eq:P1}. Finally, by \eqref{eq:T2}, \eqref{eq:T3} and \eqref{eq:NM},
\begin{align}
\nonumber E^{7/2} \left(1 + \frac{N}{M}\right)^{1/2} N^{1/2} Q^{1/2} R'^{4} \ll&\ x^{7\varepsilon/20 +\varepsilon/100} \cdot x^{1-\varepsilon} \leqslant x^{1-4\eta}, \\
\nonumber E^{3} \left(1 + \frac{N}{M}\right)^{1/2} N^{1/2 + \vartheta} Q^{1/2} R'^{4 - 3 \vartheta} \ll&\ x^{3\varepsilon/10 +\varepsilon/100} \cdot x^{1-\varepsilon} \leqslant x^{1-4\eta}.
\end{align}
This proves \eqref{eq:E1E2} and hence Theorem~\ref{thm:typeII}.
\end{proof}
\begin{proof}[Proof of Corollary~\ref{cor:symmetric}]
If $N \leqslant x^{1/2 +\varepsilon/200}$, we apply Theorem~\ref{thm:typeII} with $\varepsilon/2$ in place of $\varepsilon$: condition \eqref{eq:T1} holds, and we have
\begin{align}
\nonumber N^{1/2} Q^{1/2} R^4 \leqslant&\ x^{1/4 +\varepsilon/400}(Q R^8)^{1/2} \leqslant x^{1-\varepsilon/2}, \\
\nonumber N^{23/32} Q^{1/2} R^{107/32} \leqslant&\ x^{23/64 +\varepsilon/200}(Q^{16} R^{107})^{1/32} \leqslant x^{1-\varepsilon/2}.
\end{align}
If $N > x^{1/2 +\varepsilon/200}$, then $M \ll x^{1/2 -\varepsilon/200}$ and $M \gg x/N \geqslant Q x^{\varepsilon}$, and we apply Theorem~\ref{thm:typeII} in the same way with the roles of $M$ and $N$ reversed.
\end{proof}

\section{Proof of Theorem \ref{thm:main}}\label{sec:main}
In this section we deduce Theorem~\ref{thm:main} from Theorem~\ref{thm:typeII} and other available arithmetic information. Following \cite{MaynardLargeModuliI}, we define
$$
S_n := \mathbbm{1}_{n \equiv a (\operatorname{mod} q)} - \frac{\mathbbm{1}_{(n, q) = 1}}{\varphi(q)}, \qquad S_d(z) := \sum_{\substack{n \sim x/d \\ P^{-}(n) \geqslant z}} S_{dn}.
$$
In this and the next sections, we have $q = q_1 \cdots q_k$ in the definition of $S_n$, where $k = 2$ in this section and $k \geqslant 3$ in the next section. Let $\cQ = \cQ(\varepsilon_0)$ be the set of pairs $(Q_1, Q_2)$ of real numbers $Q_1, Q_2 \geqslant 1$ with
\begin{equation}\label{eq:region2}
Q_1^2 Q_2 \leqslant x^{1-\varepsilon_0}, \qquad Q_1^8 Q_2^{14} \leqslant x^{5-\varepsilon_0}, \qquad Q_1^{16} Q_2^{107} \leqslant x^{41/2 -\varepsilon_0}, \qquad x^{1/2 -\varepsilon} \leqslant Q_1 Q_2 \leqslant x^{8/15 -\varepsilon_0}.
\end{equation}
For $(Q_1, Q_2) \in \cQ$, we also have
\begin{align}\label{eq:conseq}
\nonumber Q_2 = \left(\frac{Q_1^8 Q_2^{14}}{(Q_1 Q_2)^8}\right)^{1/6} \leqslant x^{1/6 + 2 \varepsilon}, \qquad & Q_1^3 Q_2^2 = (Q_1^2 Q_2) (Q_1 Q_2) \leqslant x^{23/15 -2\varepsilon_0}, \\
& Q_1^{11} Q_2^{12} = (Q_1 Q_2)^{29/3} (Q_1^8 Q_2^{14})^{1/6} \leqslant x^{6 - 1/90}.
\end{align}

By a Harman's sieve argument as in [\cite{MaynardLargeModuliI}, Chapter 7], we only need to show the following four propositions:
\begin{proposition}[Type-II estimate]\label{prop:71}
Let $A > 0$ and $(Q_1, Q_2) \in \cQ$. Let $P_1, \ldots, P_J \geqslant x^{1/7 +10\varepsilon}$ be such that $P_1 \cdots P_J \ll x$ and
$$
x^{3/7 +\varepsilon} \leqslant \prod_{j \in \cJ} P_j \leqslant x^{4/7 -\varepsilon}
$$
for some subset $\cJ \subseteq \{1, \ldots, J\}$. Let $\lambda_q, \nu_q$ be complex sequences with $|\lambda_q|, |\nu_q| \leqslant \tau(q)^{B_0}$. Then we have
$$
\sum_{\substack{q_1 \sim Q_1 \\ q_2 \sim Q_2 \\ (q_1 q_2, a) = 1}} \lambda_{q_1} \nu_{q_2} \sum^{*}_{\substack{p_1, \ldots, p_J \\ p_i \sim P_i \forall i}} S_{p_1 \cdots p_J}(p_J) \ll_{A} \frac{x}{(\log x)^A}.
$$
Here $\sum^{*}$ indicates that the summation is restricted by $O(1)$ inequalities of the form $p_1^{\alpha_1} \cdots p_J^{\alpha_J} \leqslant B$. The implied constant may depend on all such exponents $\alpha_i$, but none of the quantities $B$.
\end{proposition}

\begin{proposition}[Sieve asymptotics]\label{prop:72}
Let $A > 0$ and $(Q_1, Q_2) \in \cQ$. Let $x^{3/7 +\varepsilon} \geqslant P_1 \geqslant \cdots \geqslant P_r \geqslant x^{1/7 +10\varepsilon}$ be such that $P_1 \cdots P_r \leqslant x^{5/7}$ and such that either $r = 1$ or $P_r \leqslant x^{1/4 +\varepsilon}$. Let $\lambda_q, \nu_q$ be complex sequences with $|\lambda_q|, |\nu_q| \leqslant \tau(q)^{B_0}$. Then we have
$$
\sum_{\substack{q_1 \sim Q_1 \\ q_2 \sim Q_2 \\ (q_1 q_2, a) = 1}} \lambda_{q_1} \nu_{q_2} \sum^{*}_{\substack{p_1, \ldots, p_r \\ p_i \sim P_i \forall i}} S_{p_1 \cdots p_r}(x^{1/7 +10\varepsilon}) \ll_{A} \frac{x}{(\log x)^A}.
$$
Here $\sum^{*}$ means that the summation is restricted by $O(1)$ inequalities of the form $p_1^{\alpha_1} \cdots p_r^{\alpha_r} \leqslant B$ for some constants $\alpha_1, \ldots, \alpha_r$. The implied constant may depend on all such exponents $\alpha_i$, but none of the quantities $B$. Moreover, we also have the related estimate
$$
\sum_{\substack{q_1 \sim Q_1 \\ q_2 \sim Q_2 \\ (q_1 q_2, a) = 1}} \lambda_{q_1} \nu_{q_2} S_{1}(x^{1/7 +10\varepsilon}) \ll_{A} \frac{x}{(\log x)^A}.
$$
\end{proposition}

\begin{proposition}[Numbers with four or more prime factors]\label{prop:73}
Let $A > 0$ and $(Q_1, Q_2) \in \cQ$. Let $J \geqslant 4$ and $P_1 \geqslant \cdots \geqslant P_J \geqslant x^{1/7 +10\varepsilon}$ with $P_1 \cdots P_J \asymp x$. Let $\lambda_q, \nu_q$ be complex sequences with $|\lambda_q|, |\nu_q| \leqslant \tau(q)^{B_0}$. Then we have
$$
\sum_{\substack{q_1 \sim Q_1 \\ q_2 \sim Q_2 \\ (q_1 q_2, a) = 1}} \lambda_{q_1} \nu_{q_2} \sum^{*}_{\substack{p_1, \ldots, p_J \\ p_i \sim P_i \forall i}} S_{p_1 \cdots p_J} \ll_{A} \frac{x}{(\log x)^A}.
$$
Here $\sum^{*}$ indicates that the summation is restricted by $O(1)$ inequalities of the form $p_1^{\alpha_1} \cdots p_J^{\alpha_J} \leqslant B$. The implied constant may depend on all such exponents $\alpha_i$, but none of the quantities $B$.
\end{proposition}

\begin{proposition}[Numbers with three prime factors]\label{prop:74}
Let $A > 0$ and $(Q_1, Q_2) \in \cQ$. Let $P_1, P_2, P_3 \in [x^{1/4}, x^{3/7 +\varepsilon}]$ with $P_1 P_2 P_3 \asymp x$. Let $\lambda_q, \nu_q$ be complex sequences with $|\lambda_q|, |\nu_q| \leqslant \tau(q)^{B_0}$. Then we have
$$
\sum_{\substack{q_1 \sim Q_1 \\ q_2 \sim Q_2 \\ (q_1 q_2, a) = 1}} \lambda_{q_1} \nu_{q_2} \sum^{*}_{\substack{p_1, p_2, p_3 \\ p_i \sim P_i \forall i}} S_{p_1 p_2 p_3} \ll_{A} \frac{x}{(\log x)^A}.
$$
Here $\sum^{*}$ means that the summation is restricted by $O(1)$ inequalities of the form $p_1^{\alpha_1} p_2^{\alpha_2} p_3^{\alpha_3} \leqslant B$ for some constants $\alpha_1, \alpha_2, \alpha_3$. The implied constant may depend on all such exponents $\alpha_i$, but none of the quantities $B$.
\end{proposition}

We mention that when applying Theorem~\ref{thm:typeII} and other arithmetic information proved elsewhere, we often need [\cite{MaynardLargeModuliI}, Lemma 8.11] to remove the dependencies (cross conditions) between variables. We shall omit the details of this process in this and the next sections.

\begin{proof}[Proof of Proposition~\ref{prop:71}]
We can write $S_{p_1 \cdots p_J}(p_J)$ as a sum of $S_{p_1 \cdots p_{J'}}$ with $J \leqslant J' \leqslant 6$ by a splitting argument, hence we can assume that $P_1 \cdots P_J \asymp x$. Since $Q_1 Q_2 < x^{8/15} < x^{127/224 - \varepsilon}$, [\cite{MaynardLargeModuliI}, Lemma 8.13] gives the desired result for $Q_1 Q_2 x^{\varepsilon} < \prod_{j \in \cJ} P_j \leqslant x^{4/7 -\varepsilon}$. Now we can assume that $x^{1/2} \ll \prod_{j \in \cJ} P_j \leqslant Q_1 Q_2 x^{\varepsilon}$ by symmetry. In this case we apply Theorem~\ref{thm:typeII} with $Q = Q_1$, $R = Q_2$, $m = \prod_{j \in \cJ} p_j$ and $n = \prod_{j \notin \cJ} p_j$. By \eqref{eq:region2} and \eqref{eq:conseq}, we have
\begin{align}
\nonumber & N \gg \frac{x^{1-\varepsilon}}{Q_1 Q_2} \geqslant Q_1 x^{\varepsilon_0 - \varepsilon}, \qquad N \ll x^{1/2}, \\
\nonumber & N^{1/2} Q_1^{1/2} Q_2^4 \ll x^{1/4}(Q_1^{16} Q_2^{107})^{1/32} Q_2^{21/32} \leqslant x^{1/4 + 41/64 + 7/64 -\varepsilon_0/32 + 2 \varepsilon} = x^{1 -\varepsilon_0/100}, \\
\nonumber & N^{23/32} Q_1^{1/2} Q_2^{107/32} \ll x^{23/64} (Q_1^{16} Q_2^{107})^{1/32} \leqslant x^{1 -\varepsilon_0/32}.
\end{align}
Hence \eqref{eq:T1}--\eqref{eq:T3} hold and Proposition~\ref{prop:71} follows.
\end{proof}

\begin{proof}[Proof of Proposition~\ref{prop:72}]
Note that Proposition~\ref{prop:72} is a combination of a Type-I/II estimate and a small factor Type-II estimate. We first establish an improved small factor Type-II estimate using Type-II estimates of Fouvry \cite{FouvryA1} \cite{FouvryA2} and Iwaniec--Pomyka{\l}a \cite{IwaniecPomykala}.

\begin{lemma}[Small factor Type-II estimate]\label{lem:smallfactor}
Let $A > 0$ and $Q_1, Q_2 \geqslant 1$ satisfy the first three conditions of \eqref{eq:region2}. Let $M N \asymp x$ with $x^{\varepsilon/2} \leqslant N \leqslant x^{1/6 + 11\varepsilon}$. Let $\lambda_{q_1}, \nu_{q_2}$ be $1$-bounded and let $\alpha_n, \beta_m$ be complex sequences with $|\alpha_n|, |\beta_n| \leqslant \tau(n)^{B_0}$. Suppose that $\alpha_n$ is supported on $P^{-}(n) \geqslant z_0$ and satisfies the Siegel--Walfisz condition \eqref{eq:SW}. Then we have
$$
\sum_{\substack{q_1 \sim Q_1 \\ q_2 \sim Q_2 \\ (q_1 q_2,a) = 1}} \lambda_{q_1} \nu_{q_2} \Delta(q_1 q_2) \ll_{a, \varepsilon, A, B_0} \frac{x}{(\log x)^A},
$$
where $\Delta$ is given by \eqref{eq:Delta}.
\end{lemma}
\begin{proof}[Proof of Lemma~\ref{lem:smallfactor}]
Suppose first that $x^{\varepsilon/2} \leqslant N \leqslant (Q_1 Q_2)^2 x^{-1 + 5 \varepsilon}$. We apply [\cite{FouvryA2}, Théorème(C.4)] with $R = Q_1$ and $S = Q_2$. We need
$$
\max\left(Q_1^2 Q_2, N^{5/2} Q_1 Q_2^{5/2}, N^{5/4} Q_1^{3/2} Q_2^2 \right) \leqslant x^{1 - \varepsilon/3}.
$$
By \eqref{eq:region2}, we have
\begin{align}
\nonumber & N^{5/2} Q_1 Q_2^{5/2} \leqslant (Q_1 Q_2)^{6} Q_2^{3/2} x^{-5/2 + 25\varepsilon/2} \leqslant x^{16/5 + 1/4 - 5/2 - 6 \varepsilon_0 + 16 \varepsilon} =  x^{19/20 - 6 \varepsilon_0 + 16 \varepsilon}, \\
\nonumber & N^{5/4} Q_1^{3/2} Q_2^2 \leqslant (Q_1 Q_2)^{4} Q_2^{1/2} x^{-5/4 + 25\varepsilon/4} \leqslant x^{32/15 + 1/12 - 5/4 - 4 \varepsilon_0 + 8 \varepsilon} = x^{29/30 - 4 \varepsilon_0 + 8 \varepsilon},
\end{align}
and the required condition is satisfied.

Suppose next that $(Q_1 Q_2)^2 x^{-1 + 5 \varepsilon} < N \leqslant Q_2 x^{2 \varepsilon}$. Since $Q_1^8 Q_2^{14} \leqslant x^{5 - \varepsilon_0}$, we have
$$
Q_2 x^{2 \varepsilon} < x^{5/6 - 8 \varepsilon} (Q_1 Q_2)^{-4/3}. 
$$
Now [\cite{FouvryA1}, Théorème 1] covers this case.

Suppose now that $Q_2 x^{2 \varepsilon} < N \leqslant x^{1/6 + 11\varepsilon}$. We apply [\cite{IwaniecPomykala}, Proposition] in the form of [\cite{Lichtman2}, Proposition 11.1] with $R = Q_1$ and $Q = Q_2$. We need
$$
N^6 Q_1^8 Q_2^4 < x^{5 - \varepsilon}, \qquad Q_1^2 Q_2 < N x^{1 - \varepsilon}.
$$
By \eqref{eq:region2}, we have
$$
N^6 Q_1^8 Q_2^4 = N^6 (Q_1^2 Q_2)^4 \leqslant x^{5 - 4 \varepsilon_0 + 66 \varepsilon}, \qquad Q_1^2 Q_2 \leqslant x^{1 - \varepsilon_0},
$$
and the required condition is satisfied. Now the proof of Lemma~\ref{lem:smallfactor} is completed.
\end{proof}

With Proposition~\ref{prop:71} and Lemma~\ref{lem:smallfactor}, Proposition~\ref{prop:72} can be proved in a same way of [\cite{MaynardLargeModuliI}, Proposition 7.2] with several modifications. The absolute values are replaced by bilinear weights, and the conditions on $(Q_1, Q_2)$ in \cite{MaynardLargeModuliI} are replaced by $(Q_1, Q_2) \in \cQ$. In [\cite{MaynardLargeModuliI}, Lemma 10.4] we replace the condition $N \leqslant x^{3/7 + \varepsilon} x/(Q_1 Q_2)^{15/8}$ by $N \leqslant x^{3/7 + 2 \varepsilon}$, and we have
\begin{align}
\nonumber & K L \gg \frac{x}{N} \geqslant x^{4/7 - 2 \varepsilon} \gg Q_1 Q_2 x^{\varepsilon}, \\
\nonumber & K^3 Q_1 Q_2 = \frac{K^4 Q_1 Q_2 N}{x^{2 - 2 \varepsilon}} \cdot \frac{x^{2 - 2 \varepsilon}}{N K} \ll \frac{x^{1 + 8/15 + 3/7 - \varepsilon_0 + 6 \varepsilon}}{x^{2 - 2 \varepsilon}} L x^{1 - 2 \varepsilon} < L x^{1 - 2 \varepsilon}, \\
\nonumber & K (Q_1 Q_2)^2 = \frac{K^2 (Q_1 Q_2)^2 N}{x^{2 - 2 \varepsilon}} \cdot \frac{x^{2 - 2 \varepsilon}}{N K} \ll \frac{x^{1/2 + 16/15 + 3/7 - 2 \varepsilon_0 + 4 \varepsilon}}{x^{2 - 2 \varepsilon}} L x^{1 - 2 \varepsilon} < L x^{1 - 2 \varepsilon}.
\end{align}
In [\cite{MaynardLargeModuliI}, Proposition 10.1] we replace the condition $D E P > x^{1 - 2 \varepsilon}/(Q_1 Q_2)^{15/8}$ by $D E P > x^{\varepsilon}$. Proposition~\ref{prop:72} now follows from the arguments in [\cite{MaynardLargeModuliI}, Chapter 10].
\end{proof}

\begin{proof}[Proof of Proposition~\ref{prop:73}]
Since $Q_1 Q_2 < x^{8/15} < x^{4/7 - 5\varepsilon}$, Proposition~\ref{prop:73} follows easily from the arguments in [\cite{MaynardLargeModuliI}, Chapter 9] using Proposition~\ref{prop:71} and [\cite{MaynardLargeModuliI}, Proposition 8.3], with only two modifications needed: In \textit{Case 1} ($J = 4$) we have
$$
Q_1 Q_2 K L^4 < x^{8/15} P_1 P_4 (P_2 P_3)^{3/2} \ll x^{8/15 + 1 + 3/14 + 2\varepsilon} = x^{367/210 + 2\varepsilon} < x^{57/32 - 10 \varepsilon}.
$$
In \textit{Case 2} ($J = 5$) we have
$$
Q_1 Q_2 K L^4 < x^{8/15} \frac{(P_1 P_2 P_3 P_4 P_5)^{4/3}}{(P_4 P_5)^{1/3}} \ll x^{8/15 + 4/3 - 2/21} = x^{62/35} < x^{57/32 - 10 \varepsilon}.
$$
\end{proof}

\begin{proof}[Proof of Proposition~\ref{prop:74}]
Since we have \eqref{eq:conseq}, Proposition~\ref{prop:74} follows from [\cite{Lichtman2}, Proposition 12.2], Proposition~\ref{prop:73}, a variant of [\cite{MaynardLargeModuliI}, Proposition 10.1] proved using Lemma~\ref{lem:smallfactor} (with absolute values replaced by bilinear weights, $D E P > x^{1 - 2 \varepsilon}/(Q_1 Q_2)^{15/8}$ replaced by $D E P > x^{\varepsilon}$, and conditions on $(Q_1, Q_2)$ in \cite{MaynardLargeModuliI} replaced by $(Q_1, Q_2) \in \cQ$) and the arguments in [\cite{MaynardLargeModuliI}, Chapter 11].
\end{proof}

\section{From bilinear to multilinear}\label{sec:multilinear}
In this section we show that one can get results better than Corollary~\ref{cor:multilinear} for multilinear weights by grouping divisors in different ways. We still need Harman's sieve to produce primes, but the arguments in [\cite{MaynardLargeModuliI}, Chapter 7] break down at $8/15$ since we cannot handle the sum of numbers with three balanced prime factors beyond this point. Fortunately, we can still produce primes when the Type-II range is wide enough. The following proposition shows that a Type-II range $[1/3 - \delta, 1/2 + \delta]$ for any $\delta > 0$ would be enough.

\begin{proposition}[A simple Harman's sieve]\label{prop:81}
Let $A > 0$, $k \geqslant 2$, $q_1 \cdots q_k \leqslant x^{3/5}$ and $0 < \delta \leqslant 1/100$. Let $M, N, Q_i \geqslant 1$ satisfy $M N \asymp x$. Let $\alpha_n, \beta_m, \lambda_i(q)$ be complex sequences with $|\alpha_n|, |\beta_n|, |\lambda_i(n)| \leqslant \tau(n)^{B_0}$, and assume that $\alpha_n$ is supported on integers $n$ with $P^{-}(n) \geqslant z_0$ and satisfies the Siegel--Walfisz condition \eqref{eq:SW}. Suppose that for $x^{1/3 - \delta} \leqslant N \leqslant x^{1/2 + \delta}$ we have
$$
\sum_{\substack{q_1, \ldots, q_k \\ q_i \sim Q_i \forall i \\ (q_1 \cdots q_k,a) = 1}} \lambda_1(q_1) \cdots \lambda_k(q_k) \Delta(q_1 \cdots q_k) \ll_{A} \frac{x}{(\log x)^A}
$$
with $\Delta(q)$ as in \eqref{eq:Delta}. Then,
$$
\sum_{\substack{q_1, \ldots, q_k \\ q_i \sim Q_i \forall i \\ (q_1 \cdots q_k,a) = 1}} \lambda_1(q_1) \cdots \lambda_k(q_k) \sum_{p \sim x} S_p \ll_{A} \frac{x}{(\log x)^A}.
$$
\end{proposition}
\begin{proof}
Since $\sum_{p \sim x} S_p = S_1(2x^{1/2})$, we can apply Buchstab's identity repeatedly (see [\cite{LRB2602}, Lemma 2.13]) to get two types of sums:
$$
\sum_{\substack{x^{\delta/4} \leqslant p_j < \cdots < p_1 < x^{1/2} \\ p_1 \cdots p_j < 2 x^{1/3 - \delta}}} S_{p_1 \cdots p_j}(x^{\delta/4}) \qquad \text{and} \qquad \sum_{\substack{x^{\delta/4} \leqslant p_{j+1} < \cdots < p_1 < x^{1/2} \\ p_1 \cdots p_j < 2 x^{1/3 - \delta} \leqslant p_1 \cdots p_{j+1}}} S_{p_1 \cdots p_{j+1}}(p_{j+1}).
$$
For the second type of sums, we have
$$
2 x^{1/3 - \delta} \leqslant p_1 \cdots p_{j+1} < p_1 \cdots p_j^2 \leqslant (p_1 \cdots p_j)^2 < 4 x^{2/3 - 2 \delta}.
$$
Since the sum counts numbers of the form $p_1 \cdots p_{j+1} d$ with $P^{-}(d) \geqslant p_{j+1}$, we can use our assumption with $n = p_1 \cdots p_{j+1}$ when $p_1 \cdots p_{j+1} \leqslant 2 x^{1/2}$ and $n = d$ when $p_1 \cdots p_{j+1} > 2 x^{1/2}$. Clearly $\alpha_n$ satisfies the Siegel--Walfisz condition \eqref{eq:SW} in both case. Now we have
$$
\sum_{\substack{x^{\delta/4} \leqslant p_{j+1} < \cdots < p_1 < x^{1/2} \\ p_1 \cdots p_j < 2 x^{1/3 - \delta} \leqslant p_1 \cdots p_{j+1}}} S_{p_1 \cdots p_{j+1}}(p_{j+1}) \ll \frac{x}{(\log x)^A}
$$
for all $j$.

For the first type of sums, we apply [\cite{MaynardLargeModuliI}, Lemma 10.2] with $z_1 = x^{\delta/4}$, $z_2 = z_0$ and $y = 2 x^{1/3 - \delta}$ to get another two types of sums:
$$
\sum_{\substack{x^{\delta/4} \leqslant p_j < \cdots < p_1 < x^{1/2} \\ p_1 \cdots p_j < 2 x^{1/3 - \delta}}} \sum_{d < 2 x^{1/3 - \delta}} \alpha_d \sum_{\substack{p_1 \cdots p_j d m \sim x \\ P^{-}(m) \geqslant z_0}} S_{p_1 \cdots p_j d m}
$$
and
$$
\sum_{\substack{x^{\delta/4} \leqslant p_j < \cdots < p_1 < x^{1/2} \\ p_1 \cdots p_j < 2 x^{1/3 - \delta}}} \sum_{z_0 \leqslant p < x^{\delta/4}} \sum_{\substack{d < 2 x^{1/3 - \delta} \leqslant d p \\ P^{-}(d) > p }} \beta_d \sum_{\substack{p_1 \cdots p_j p d m \sim x \\ P^{-}(m) \geqslant p}} S_{p_1 \cdots p_j p d m},
$$
where $\alpha_d$, $\beta_d$ are $1$-bounded and supported on $P^{-}(d) \geqslant z_0$. For the sum over $p_1 \cdots p_j p d m$, we have
$$
2 x^{1/3 - \delta} \leqslant d p < 2 x^{1/3 - \delta + \delta/4}.
$$
Similarly, we can use our assumption with $n = d p$ and get
$$
\sum_{\substack{x^{\delta/4} \leqslant p_j < \cdots < p_1 < x^{1/2} \\ p_1 \cdots p_j < 2 x^{1/3 - \delta}}} \sum_{z_0 \leqslant p < x^{\delta/4}} \sum_{\substack{d < 2 x^{1/3 - \delta} \leqslant d p \\ P^{-}(d) > p }} \beta_d \sum_{\substack{p_1 \cdots p_j p d m \sim x \\ P^{-}(m) \geqslant p}} S_{p_1 \cdots p_j p d m} \ll \frac{x}{(\log x)^A}
$$
for all $j$. The same bound also holds for the part of the sum over $p_1 \cdots p_j d m$ with $p_1 \cdots p_j d \geqslant 2 x^{1/3 - \delta}$, since we have
$$
2 x^{1/3 - \delta} \leqslant p_1 \cdots p_j d < 4 x^{2/3 - 2 \delta}.
$$
Finally, we can use the Fundamental Lemma of the sieve [\cite{MaynardLargeModuliI}, Lemma 10.3] to handle the part of the sum over $p_1 \cdots p_j d m$ with $p_1 \cdots p_j d < 2 x^{1/3 - \delta}$. We have
\begin{align}
\nonumber &\ \sum_{\substack{x^{\delta/4} \leqslant p_j < \cdots < p_1 < x^{1/2} \\ p_1 \cdots p_j d < 2 x^{1/3 - \delta}}} \alpha_d \sum_{\substack{p_1 \cdots p_j d m \sim x \\ P^{-}(m) \geqslant z_0}} S_{p_1 \cdots p_j d m} \\
\nonumber =&\ \sum_{l < 2 x^{1/3 - \delta}} c_l \sum_{\substack{l m \sim x \\ P^{-}(m) \geqslant z_0}} S_{l m}, \qquad c_l = \sum_{\substack{l = p_1 \cdots p_j d \\ x^{\delta/4} \leqslant p_j < \cdots < p_1 < x^{1/2}}} \alpha_d.
\end{align}
Using [\cite{MaynardLargeModuliI}, Lemma 10.3] for each $\sum_{\substack{l m \sim x \\ P^{-}(m) \geqslant z_0}} S_{l m}$ and summing over $q_1, \ldots, q_k$ and $l$, we obtain
\begin{align}
\nonumber \sum_{\substack{q_1, \ldots, q_k \\ q_i \sim Q_i \forall i \\ (q_1 \cdots q_k,a) = 1}} \lambda_1(q_1) \cdots \lambda_k(q_k) & \sum_{l < 2 x^{1/3 - \delta}} c_l \sum_{\substack{l m \sim x \\ P^{-}(m) \geqslant z_0}} S_{l m} \\
\nonumber \ll&\ x (\log x)^{O(1) - \log \log x} + x^{3/5 + 1/3} y_0 (\log x)^{O(1)} \\ 
\nonumber \ll&\ \frac{x}{(\log x)^A},
\end{align}
since $q_1 \cdots q_k \leqslant x^{3/5}$. Now the proof of Proposition~\ref{prop:81} is completed.
\end{proof}

With Proposition~\ref{prop:81}, we reduce the original problem to proving a Type-II range $[1/3 - \delta, 1/2 + \delta]$. Apart from Theorem~\ref{thm:typeII}, we also need the following Type-II information input come from the work of Lichtman \cite{Lichtman3} and Pascadi \cite{Pascadi2}.
\begin{lemma}[Type-II estimate for trilinear weights]\label{lem:typeIImulti}
Let $\vartheta$ be admissible. Let $a$ be a fixed non-zero integer and let $\varepsilon > 0$, $A > 0$, $B_0 > 0$. Let $M, N, Q, R, S \geqslant 1$ satisfy $M N \asymp x$ and any of the following two sets of conditions holds:
\begin{align}
\nonumber & Q x^{\varepsilon} \leqslant N, \qquad N^2 R S^2 \leqslant N^{1 - \varepsilon}, \qquad N Q R^5 S^2 \leqslant N^{2 - \varepsilon}, \qquad N^{1 + 2 \vartheta} Q R^{5 - 2 \vartheta} S^{2 + 2 \vartheta} \leqslant N^{2 - \varepsilon}; \\
\nonumber & Q x^{\varepsilon} \leqslant N, \qquad N^2 R S^2 \leqslant N^{1 - \varepsilon}, \qquad N^{\frac{2 + \vartheta}{2 - 2 \vartheta}} R^{\frac{4 - 5 \vartheta}{2 - 2 \vartheta}} S \leqslant N^{1 - \varepsilon}, \qquad N^{\frac{1 + \vartheta}{1 - \vartheta}} Q^{\frac{1 - 3 \vartheta}{1 - \vartheta}} R^5 S^2 \leqslant N^{2 - \varepsilon}.
\end{align}
Let $\alpha_n, \beta_m, \lambda_q, \nu_r, \eta_s$ be complex sequences with $|\alpha_n|, |\beta_n|, |\lambda_n|, |\nu_n|, |\eta_n| \leqslant \tau(n)^{B_0}$, and assume that $\alpha_n$ is supported on integers $n$ with $P^{-}(n) \geqslant z_0$ and satisfies the Siegel--Walfisz condition \eqref{eq:SW}. Then, with $\Delta(q)$ as in \eqref{eq:Delta},
$$
\sum_{q \sim Q} \sum_{r \sim R} \sum_{\substack{s \sim S \\ (q r s,a) = 1}} \lambda_q \nu_r \eta_s \Delta(q r s) \ll_{a, \varepsilon, A, B_0} \frac{x}{(\log x)^A}.
$$
Moreover, suppose that any of the following two sets of conditions holds:
\begin{align}
\nonumber & Q x^{\varepsilon} \leqslant N, \qquad N^2 R \leqslant N^{1 - \varepsilon}, \qquad N Q R^5 \leqslant N^{2 - \varepsilon}, \qquad N^{1 + 2 \vartheta} Q R^{5 - 2 \vartheta} \leqslant N^{2 - \varepsilon}; \\
\nonumber & Q x^{\varepsilon} \leqslant N, \qquad N^2 R \leqslant N^{1 - \varepsilon}, \qquad N^{\frac{2 + \vartheta}{2 - 2 \vartheta}} R^{\frac{4 - 5 \vartheta}{2 - 2 \vartheta}} \leqslant N^{1 - \varepsilon}, \qquad N^{\frac{1 + \vartheta}{1 - \vartheta}} Q^{\frac{1 - 3 \vartheta}{1 - \vartheta}} R^5 \leqslant N^{2 - \varepsilon}.
\end{align}
Then, with $\Delta(q)$ as in \eqref{eq:Delta},
$$
\sum_{q \sim Q} \sum_{\substack{r \sim R \\ (q r,a) = 1}} \lambda_q \nu_r \Delta(q r) \ll_{a, \varepsilon, A, B_0} \frac{x}{(\log x)^A}.
$$
\end{lemma}
\begin{proof}
This follows from [\cite{Lichtman3}, Proposition 5.2] and [\cite{Pascadi2}, Proposition 4.4].
\end{proof}

Now, we have three different Type-II ranges generated by Theorem~\ref{thm:typeII} and Lemma~\ref{lem:typeIImulti}:
\begin{align}
& \left[Q x^{\varepsilon},\ x^{-\varepsilon} \min\left(x^{1/2 + 101 \varepsilon/100},\ x^2 Q^{-1} R^{-8},\ x^{\frac{2}{1 + 2 \vartheta}} Q^{-\frac{1}{1 + 2 \vartheta}} R^{-\frac{8 - 6 \vartheta}{1 + 2 \vartheta}} \right) \right], \label{eq:TypeIIrange14} \\
& \left[Q x^{\varepsilon},\ x^{-\varepsilon} \min\left(x^{1/2} R^{-1/2} S^{-1},\ x^2 Q^{-1} R^{-5} S^{-2},\ x^{\frac{2}{1 + 2 \vartheta}} Q^{-\frac{1}{1 + 2 \vartheta}} R^{-\frac{5 - 2 \vartheta}{1 + 2 \vartheta}} S^{-\frac{2 + 2 \vartheta}{1 + 2 \vartheta}} \right) \right], \label{eq:TypeIIrangeLichtman} \\
& \left[Q x^{\varepsilon},\ x^{-\varepsilon} \min\left(x^{1/2} R^{-1/2} S^{-1},\ x^{\frac{2 - 2 \vartheta}{2 + \vartheta}} R^{\frac{4 - 5 \vartheta}{2 + \vartheta}} S^{-\frac{2 - 2 \vartheta}{2 + \vartheta}},\ x^{\frac{2 - 2 \vartheta}{1 + \vartheta}} Q^{-\frac{1 - 3 \vartheta}{1 + \vartheta}} R^{-\frac{5 - 5 \vartheta}{1 + \vartheta}} S^{-\frac{2 - 2 \vartheta}{1 + \vartheta}} \right) \right]. \label{eq:TypeIIrangePascadi}
\end{align}
Hence, we can get multilinear results for primes if we can group $Q_1, \ldots, Q_k$ into $Q, R$ or $Q, R, S$ in different ways such that the union of all available Type-II ranges from \eqref{eq:TypeIIrange14}--\eqref{eq:TypeIIrangePascadi} covers $[1/3 - \delta, 1/2 + \delta]$ for some $\delta > 0$. This leads to a set of linear inequalities on $\theta_i := \log Q_i / \log x$. The complete explicit region of the solution of such inequalities is extremely complicated, see [\cite{LRB2602}, Section 5] for similar things; we shall prove a subregion of all solutions for the trilinear case and only indicate special solution values for the quadrilinear and quintilinear cases. Note that we do not fully optimize the parameters, and further improvements on the exponents $145/266$, $457/814$ and $3614/6353$ are still possible.

\begin{theorem}\label{thm:trilinear}
Let $a$ be a fixed non-zero integer and let $\varepsilon > 0$, $A > 0$, $B_0 > 0$. Let $Q_1, Q_2, Q_3 \geqslant 1$ satisfy
\begin{align}
& Q_2 Q_3 \leqslant x^{1/3 - \varepsilon}, \label{eq:tricond1} \\
& Q_1^3 Q_3^2 \leqslant x^{1 - \varepsilon}, \qquad Q_1^{96} Q_2^{16} Q_3^{39} \leqslant x^{32 - \varepsilon}, \label{eq:tricond2} \\
& Q_1^2 Q_2^3 \leqslant x^{1 - \varepsilon}, \qquad Q_1^{39} Q_2^{96} Q_3^{16} \leqslant x^{32 - \varepsilon}, \label{eq:tricond3} \\
& Q_1 Q_2 Q_3^8 \leqslant x^{3/2 - \varepsilon}, \qquad Q_1^{16} Q_2^{16} Q_3^{107} \leqslant x^{41/2 - \varepsilon}. \label{eq:tricond4}
\end{align}
Let $\lambda_1, \lambda_2, \lambda_3$ be complex sequences with $|\lambda_i(q)| \leqslant \tau(q)^{B_0}$. Then
$$
\sum_{\substack{q_1 \leqslant Q_1 \\ q_2 \leqslant Q_2 \\ q_3 \leqslant Q_3 \\ (q_1 q_2 q_3, a) = 1}} \lambda_1(q_1) \lambda_2(q_2) \lambda_3(q_3) \left( \pi(x; q_1 q_2 q_3, a) - \frac{\pi(x)}{\varphi(q_1 q_2 q_3)} \right) \ll_{a, \varepsilon, A, B_0} \frac{x}{(\log x)^A}.
$$
\end{theorem}
\begin{remark*}
Under the conditions \eqref{eq:tricond1}--\eqref{eq:tricond4}, we have
$$
(Q_1 Q_2 Q_3)^{133} = (Q_1^3 Q_3^2)^{13} (Q_1^2 Q_2^3)^{39} Q_1^{16} Q_2^{16} Q_3^{107} \leqslant x^{145/2 - 53 \varepsilon},
$$
hence $Q_1 Q_2 Q_3 \leqslant x^{145/266 - \varepsilon}$, and the maximum value is attained at
$$
Q_1 = x^{230/931 - \varepsilon}, \qquad Q_2 = x^{157/931 - \varepsilon}, \qquad Q_3 = x^{241/1862 - \varepsilon}.
$$
\end{remark*}
\begin{proof}
Let $\varepsilon' = \varepsilon/1000$. We shall use three different Type-II ranges to cover $[1/3 - \delta, 1/2 + \delta]$: The first range is \eqref{eq:TypeIIrangeLichtman} with $Q = Q_2 Q_3$, $R = Q_1$ and $S = 1$, hence \eqref{eq:tricond1} is needed to cover $1/3 - \delta$. The second range is \eqref{eq:TypeIIrangeLichtman} with $Q = Q_1 Q_3$, $R = Q_2$ and $S = 1$, and we need
$$
Q_1 Q_3 x^{2 \varepsilon'} \leqslant \min\left(x^{1/2} Q_1^{-1/2},\ x^2 (Q_2 Q_3)^{-1} Q_1^{-5},\ x^{32/23} (Q_2 Q_3)^{-16/23} Q_1^{-73/23} \right).
$$
This is satisfied by \eqref{eq:tricond2}, since we have
$$
Q_1^6 Q_2 Q_3^2 \ll Q_1^6 Q_2 Q_3^{39/16} = (Q_1^{96} Q_2^{16} Q_3^{39})^{1/16} \leqslant x^{2 - \varepsilon/16}.
$$
The third range is \eqref{eq:TypeIIrange14} with $Q = Q_1 Q_2$ and $R = Q_3$, and we need
$$
Q_1 Q_2 x^{2 \varepsilon'} \leqslant \min\left(x^{1/2} Q_2^{-1/2},\ x^2 (Q_1 Q_3)^{-1} Q_2^{-5},\ x^{32/23} (Q_1 Q_3)^{-16/23} Q_2^{-73/23} \right).
$$
This is satisfied by \eqref{eq:tricond3}, since we have
$$
Q_1^2 Q_2^6 Q_3 \ll Q_1^{39/16} Q_2^6 Q_3 = (Q_1^{39} Q_2^{96} Q_3^{16})^{1/16} \leqslant x^{2 - \varepsilon/16}.
$$
We also need
$$
x^{1/2 + \varepsilon'} \leqslant \min\left(x^2 (Q_1 Q_2)^{-1} Q_3^{-8},\ x^{32/23} (Q_1 Q_2)^{-16/23} Q_3^{-107/23} \right)
$$
in order to cover $1/2 + \delta$, and this is satisfied by \eqref{eq:tricond4}.
\end{proof}

\begin{theorem}\label{thm:345linear}
Let $a$ be a fixed non-zero integer and let $\varepsilon > 0$, $A > 0$, $B_0 > 0$. Let $k \geqslant 3$ and $Q_1, \ldots, Q_k \geqslant 1$ satisfy any of the following three conditions:
\begin{align}
k = 3, \qquad & Q_1 = x^{83/250 - \varepsilon}, \quad Q_2 = x^{69/500 - \varepsilon}, \quad Q_3 = x^{13/200 - \varepsilon}; \label{eq:3linearcond} \\
k = 4, \qquad & Q_1 = x^{439/1628 - \varepsilon}, \quad Q_2 = x^{225/1628 - \varepsilon}, \quad Q_3 = x^{175/1628 - \varepsilon}, \quad Q_4 = x^{75/1628 - \varepsilon}; \label{eq:4linearcond} \\
\nonumber k = 5, \qquad & Q_1 = x^{1789/6353 - \varepsilon}, \quad Q_2 = x^{750/6353 - \varepsilon}, \quad Q_3 = x^{675/6353 - \varepsilon}, \\
& Q_4 = x^{225/6353 - \varepsilon}, \quad Q_5 = x^{175/6353 - \varepsilon}. \label{eq:5linearcond}
\end{align}
Let $\lambda_1, \ldots, \lambda_k$ be complex sequences with $|\lambda_i(q)| \leqslant \tau(q)^{B_0}$. Then
$$
\sum_{\substack{q_1, \ldots, q_k \\ q_i \leqslant Q_i \forall i \\ (q_1 \cdots q_k,a) = 1}} \lambda_1(q_1) \cdots \lambda_k(q_k) \left( \pi(x; q_1 \cdots q_k, a) - \frac{\pi(x)}{\varphi(q_1 \cdots q_k)} \right) \ll_{a, \varepsilon, A, B_0} \frac{x}{(\log x)^A}.
$$
\end{theorem}
\begin{remark*}
One can easily see that the first case of Theorem~\ref{thm:345linear} is not covered in Theorem~\ref{thm:trilinear}. If it is covered, then we must have $Q_1 = x^{83/250 - \varepsilon}$ in Theorem~\ref{thm:trilinear} in order to satisfy \eqref{eq:tricond1}. But now we have
$$
Q_1^3 Q_3^2 \geqslant x^{249/250 + 13/100 - 5 \varepsilon} = x^{563/500 - 5 \varepsilon} > x,
$$
contradicting with \eqref{eq:tricond2}.
\end{remark*}
\begin{proof}
We group $Q_1, \ldots, Q_k$ in different ways and apply Type-II ranges \eqref{eq:TypeIIrange14} and \eqref{eq:TypeIIrangeLichtman} to cover $[1/3 - \delta, 1/2 + \delta]$. One can see Table~\ref{tab:ranges} for the numerical boundaries of the ranges used.
\end{proof}
\begin{table}[ht]
\renewcommand{\arraystretch}{1.25}
\begin{tabular}{llll}
\hline
$k$ & Type-II ranges used & Lower bounds & Upper bounds \\
\hline
$3$ & \eqref{eq:TypeIIrangeLichtman}, $Q = Q_1$, $R = Q_2 Q_3$, $S = 1$ & $83/250$ & $797/2000$ \\
& \eqref{eq:TypeIIrange14}, $Q = Q_1 Q_3$, $R = Q_2$ & $397/1000$ & $5441/11500$ \\
& \eqref{eq:TypeIIrange14}, $Q = Q_1 Q_2$, $R = Q_3$ & $47/100$ & $1/2$ \\
\hline
$4$ & \eqref{eq:TypeIIrangeLichtman}, $Q = Q_1 Q_4$, $R = Q_2 Q_3$, $S = 1$ & $257/814$ & $307/814$ \\
& \eqref{eq:TypeIIrangeLichtman}, $Q = Q_1 Q_3$, $R = Q_2 Q_4$, $S = 1$ & $307/814$ & $166/407$ \\
& \eqref{eq:TypeIIrangeLichtman}, $Q = Q_1 Q_2$, $R = Q_3 Q_4$, $S = 1$ & $166/407$ & $689/1628$ \\
& \eqref{eq:TypeIIrange14}, $Q = Q_1 Q_3 Q_4$, $R = Q_2$ & $689/1628$ & $739/1628$ \\
& \eqref{eq:TypeIIrange14}, $Q = Q_1 Q_2 Q_4$, $R = Q_3$ & $739/1628$ & $1/2$ \\
\hline
$5$ & \eqref{eq:TypeIIrangeLichtman}, $Q = Q_1 Q_4$, $R = Q_2 Q_3 Q_5$, $S = 1$ & $2014/6353$ & $54272/146119$ \\
& \eqref{eq:TypeIIrangeLichtman}, $Q = Q_1 Q_4 Q_5$, $R = Q_2 Q_3$, $S = 1$ & $2189/6353$ & $2464/6353$ \\
& \eqref{eq:TypeIIrangeLichtman}, $Q = Q_1 Q_3$, $R = Q_2 Q_4 Q_5$, $S = 1$ & $2464/6353$ & $5203/12706$ \\
& \eqref{eq:TypeIIrangeLichtman}, $Q = Q_1 Q_2$, $R = Q_3 Q_4 Q_5$, $S = 1$ & $2539/6353$ & $2639/6353$ \\
& \eqref{eq:TypeIIrangeLichtman}, $Q = Q_1 Q_3 Q_5$, $R = Q_2 Q_4$, $S = 1$ & $2639/6353$ & $2689/6353$ \\
& \eqref{eq:TypeIIrangeLichtman}, $Q = Q_1 Q_3 Q_4$, $R = Q_2 Q_5$, $S = 1$ & $2689/6353$ & $2714/6353$ \\
& \eqref{eq:TypeIIrange14}, $Q = Q_1 Q_2 Q_5$, $R = Q_3 Q_4$ & $2714/6353$ & $2764/6353$ \\
& \eqref{eq:TypeIIrange14}, $Q = Q_1 Q_2 Q_4$, $R = Q_3 Q_5$ & $2764/6353$ & $68122/146119$ \\
& \eqref{eq:TypeIIrange14}, $Q = Q_1 Q_2 Q_4 Q_5$, $R = Q_3$ & $2939/6353$ & $1/2$ \\
\hline
\end{tabular}
\medskip
\caption{\small{Type-II ranges used in the proof of Theorem~\ref{thm:345linear}}}\label{tab:ranges}
\end{table}

Finally, we show that the maximum size of $Q_1 \cdots Q_k$ can be arbitrarily close to $x^{4/7}$ when $k$ is large enough. Assume that $\theta_1 + \cdots + \theta_k = \theta$, where $1/2 \leqslant \theta < 4/7$. Let $S = 1$ and write $Q = x^{\theta'}$, $R = x^{\theta''}$. Then, the Type-II ranges \eqref{eq:TypeIIrange14} and \eqref{eq:TypeIIrangeLichtman} becomes
\begin{equation}\label{eq:typeII14-47}
\left[\theta' + \varepsilon,\ \min\left(\frac{1}{2} + \frac{101\varepsilon}{100},\ 2 - \theta' - 8 \theta'',\ \frac{2 - \theta' - (8 - 6 \vartheta) \theta''}{1 + 2 \vartheta} \right) - \varepsilon \right]
\end{equation}
and
\begin{equation}\label{eq:typeIILichtman-47}
\left[\theta' + \varepsilon,\ \min\left(\frac{1}{2} - \frac{\theta''}{2},\ 2 - \theta' - 5 \theta'',\ 1 - \frac{\theta'}{2} - 2 \theta'' \right) - \varepsilon \right]. \qquad (\vartheta = 1/2)
\end{equation}
Let $\nu \in [1/3 - \delta, 1/2 + \delta]$ and put
$$
\nu_0 := \max\left(\frac{4 \theta - 1}{3},\ \frac{(4 - 3 \vartheta) \theta - 1}{3 - 4 \vartheta} \right).
$$
Since $\theta < 4/7$, we have $\nu_0 < 1 - \theta$. Let $\theta_1 = \cdots = \theta_k = h := \theta/k$, and let $k$ be large such that $h < \min(1/1000, (1 - \theta - \nu_0)/3)$. Let $m$ be the unique integer such that $m h < \nu < (m + 1)h$, and let $Q = Q_1 \cdots Q_m$ and $R = Q_{m + 1} \cdots Q_k$. Thus, we have $\nu - h \leqslant \theta' < \nu$ and $\theta'' \leqslant \theta - \nu + h$.

If $\nu < 1 - \theta - h$, then we have
$$
\frac{1 - \theta''}{2} \geqslant \frac{1 - \theta + v - h}{2} > v,
$$
$$
2 - \theta' - 5 \theta'' \geqslant 2 - 5 \theta + 4 (\nu - h) > \nu \qquad \text{and} \qquad 1 - \frac{\theta'}{2} - 2 \theta'' \geqslant 1 - 2 \theta + \frac{3(\nu - h)}{2} > \nu
$$
since $\nu \geqslant 1/3 - \delta$, $\theta < 4/7$ and $h < 1/1000$. Now $\nu$ is in \eqref{eq:typeIILichtman-47}.

If $\nu > \nu_0 + 2 h$, then we have
$$
2 - \theta' - 8 \theta'' \geqslant 2 - 8 \theta + 7 (\nu - h) > \nu
$$
and
$$
2 - \theta' - (8 - 6 \vartheta) \theta'' \geqslant 2 - (8 - 6 \vartheta) \theta + (7 - 6 \vartheta) (\nu - h) > (1 + 2 \vartheta) \nu.
$$
Now $\nu$ is in \eqref{eq:typeII14-47}. Since $\nu_0 + 2 h < 1 - \theta - h$, every $\nu \in [1/3 - \delta, 1/2 + \delta]$ is in \eqref{eq:typeII14-47} or \eqref{eq:typeIILichtman-47}, and the following Theorem is proved.

\begin{theorem}\label{thm:approach47}
Let $a$ be a fixed non-zero integer and let $\varepsilon > 0$, $A > 0$, $B_0 > 0$. Let $\lambda_i$ be complex sequences with $|\lambda_i(q)| \leqslant \tau(q)^{B_0}$. Then for every $1/2 \leqslant \theta < 4/7$, there exists a positive integer $k$ such that
$$
\sum_{\substack{q_1, \ldots, q_k \\ q_i \leqslant x^{\theta/k} \forall i \\ (q_1 \cdots q_k,a) = 1}} \lambda_1(q_1) \cdots \lambda_k(q_k) \left( \pi(x; q_1 \cdots q_k, a) - \frac{\pi(x)}{\varphi(q_1 \cdots q_k)} \right) \ll_{a, \varepsilon, A, B_0} \frac{x}{(\log x)^A}.
$$
\end{theorem}

\section*{AI declaration}
Claude Fable 5.2 (Anthropic) was used in the proof of parts of Theorem~\ref{thm:typeII} and in the preparation of the manuscript, including language editing and formatting. The author would also like to thank Shuishan Jiang for his help on the usage of AI.

\bibliographystyle{plain}
\bibliography{bib}
\end{document}